\documentclass[10pt,twoside]{article}

\usepackage[centertags]{amsmath}
\usepackage{amsfonts}
\usepackage{amssymb}
\usepackage{amsthm}
\usepackage{epsfig}
\usepackage{color}
\usepackage{bbm}
\usepackage{natbib}
\allowdisplaybreaks[4]
\date{}

\definecolor{c20}{rgb}{0.,0.7,0.}
\definecolor{c30}{rgb}{0.,0.,1.}
\definecolor{c40}{rgb}{1,0.1,0.7}
\definecolor{c50}{rgb}{1,0,0}
\definecolor{c60}{rgb}{1,0.9,0.1}

\begin{document}
\baselineskip 15pt \setcounter{page}{1}
\title{\bf \Large
On the maxima of nonstationary random fields subject to missing observations
\thanks{Research supported by Innovation of Jiaxing City: a program to support the talented persons, National Bureau of Statistics of China (No. 2020LY031)
and Project of new economy research center of Jiaxing City (No. WYZB202254, WYZB202257).}}
\author{{\small Shengchao Zheng,\ \ Zhongquan Tan}\footnote{ E-mail address:  tzq728@163.com }\\
\\
{\small\it  College of Data Science, Jiaxing University, Jiaxing 314001, PR China}\\
}
 \maketitle
 \baselineskip 15pt

\begin{quote}
{\bf Abstract:}\ \ Motivated by the papers of \cite{Mladenovc_Piterbarg2006}, \cite{Krajka2011} and \cite{Pereira_Tan2017},  we study the limit properties for the maxima from nonstationary random fields subject to missing observations and obtain the weakly convergence and almost sure convergence results for these maxima. Some examples such as Gaussian random fields, $chi$-random fields and Gaussian order statistics fields are given to illustrate the obtained results.

{\bf Key Words:}\ extreme value theory; nonstationary random fields; missing observations; almost sure central limit theorem

{\bf AMS Classification:}\ \ Primary 60G70; secondary 60G60

\end{quote}

\section{Introduction}
Suppose that $\{X_n,n\geq 1\}$ \ is a sequence of stationary random variables with the marginal distribution function $F(x)$ and satisfies the weak dependence condition $D(u_{n})$ and the local dependence condition $D'(u_{n})$ (see e.g., \cite{Leadbetter_Lindgren_Root1983} for the definitions), where $u_{n}=u_{n}(x)=a_{n}^{-1}x+b_{n}$ is a sequence of constants with $a_{n}>0$ and $b_{n}\in \mathbb{R}$. If $u_{n}(x)$ satisfies $n(1-F(u_{n}(x)))\rightarrow -\log G(x)$, as $n\rightarrow\infty$, then for any $x\in \mathbb{R}$, we have
 \begin{eqnarray}
\label{01}
\lim_{n\rightarrow\infty}P\left(M_{n}\leq u_{n}(x)\right)= G(x),
\end{eqnarray}
where $M_{n}=\max\{X_{k}, k=1,2,\ldots,n\}$, and $G(x)$ is of one of the three types of extreme value distributions. The more details for (\ref{01}) can be found in the monographs \citep{Leadbetter_Lindgren_Root1983,Piterbarg_1996}.
Missing observations may occur randomly and cause serious consequences in practical applications. It is very important to study the impact induced by the missing observations on the maxima in extreme value theory.  Suppose $\varepsilon_{k}$ is the indicator of the event whether the random variable $X_{k}$ is observed (in other word, $\varepsilon_{k}$ is a Bernoulli random variable), and let $S_{n}=\sum_{k\leq n}\varepsilon_{k}$. Furthermore, we  assume that
$\{\varepsilon_{n},{n\geq 1}\}$ are independent of $\{X_{n},{n\geq 1}\}$
with $S_{n}$ satisfying
   \begin{eqnarray*}
\frac{S_{n}}{n}\stackrel{P}{\longrightarrow} \lambda,\ \ \mbox{as}\quad n\rightarrow\infty,
\end{eqnarray*}
where $\lambda$ is a random or nonrandom variable.

 For  a constant  $\lambda\in (0,1)$, \cite{Mladenovc_Piterbarg2006} studied the asymptotic distribution of maximum from stationary sequences and its maximum when the sequence is subject to random missing observations under conditions $D(u_{n},v_{n})$ and $D'(u_{n})$, and derived the following result for any $x < y \in \mathbb{R}$,
\begin{eqnarray}
\label{03}
\lim_{n\rightarrow\infty}P\left( \widetilde{M}_{n}\leq u_{n}(x), M_{n}\leq v_{n}(y)\right)=G^{\lambda}(x)G^{1-\lambda}(y),
\end{eqnarray}
where $\widetilde{M}_{n}=\max\{X_{k}, \varepsilon_{k}=1, k=1,2,\ldots,n\}$ denotes the  maximum subject to random missing observations.
\cite{Cao_Peng2011} and \cite{Peng_Cao_Nadarajah2010} extended the result  (\ref{03}) to Gaussian cases. The generalization of (\ref{03}) to autoregressive process and linear process can be found in \cite{Glava_Mladenov2017} and \cite{Glava_Mladenovi2020}, respectively, and to nonstationary random fields can be found in \cite{Panga_Pereira2018}.
\cite{Tong_Peng2011} also studied
the almost sure limit theorem for these maxima and got
\begin{equation}
\lim_{n\rightarrow\infty}\frac{1}{\log n}\sum_{k=1}^{n}\frac{1}{k}\mathbf{\mathbbm{1}}_{\left\{\widetilde{M}_{n}\leq u_{n}(x), M_{n}\leq v_{n}(y)\right\}}= G^{\lambda}(x)G^{1-\lambda}(y) \ \ \ \   a.s.,
\end{equation}
for any $x<y\in \mathbb{R}$.

When $\lambda$  is  a random variable, \cite{Krajka2011} obtained the following result for any $x< y \in \mathbb{R}$
\begin{eqnarray}
\label{04}
\lim_{n\rightarrow\infty}P\left( \widetilde{M}_{n}\leq u_{n}(x), M_{n}\leq v_{n}(y)\right)= E[G^{\lambda}(x)G^{1-\lambda}(y)].
\end{eqnarray}
\cite{Hashorva_Peng_Weng2013} extended the results of (\ref{04}) to weakly and strongly dependent Gaussian sequences.

To our knowledge, the studies on the maxima subject to missing observations mainly focus on the case when $\lambda$ is a constant.
In this paper, we will extend the above limit results on the maxima subject to missing observations to nonstationary random fields when
 $\lambda$ is a random variable.

We first generalize the dependent condition $D(u_{n})$ and the local dependence condition $D'(u_{n})$ for nonstationary random fields.
These conditions strongly rely on the random variables, the partition of index set and the levels $u_{n}$.

For sake of simplicity we only deal with the
two-dimensional case and we claim that the results for higher dimensions can be derived by similar discussions. We first introduce some notations and notions used in the paper. For $\mathbf{i}=(i_1,i_2)$
 and $\mathbf{j}=(j_1,j_2)$, $\mathbf{i} \leq \mathbf{j}$ means $i_k \leq j_k, k=1,2$ and $\mathbf{n}=
 (n_1,n_2)\rightarrow \infty$ means $n_k\rightarrow \infty, k=1,2 $, respectively.   Let $\mathbf{X}=\{X_{\mathbf{n}},\mathbf{n}\geq \mathbf{1}\}$ be a sequence of  nonstationary random fields and $\{u_{\mathbf{n,i}}, \mathbf{i\leq n}\}_{\mathbf{n \geq 1}}$ and $\{v_{\mathbf{n,i}},\mathbf{i\leq n}\}_{\mathbf{n \geq 1}}$ be two sequences of real numbers satisfying $v_{\mathbf{n,i}}\leq u_{\mathbf{n,i}}$ for all $\mathbf{i\leq n}$.
Suppose that $P(X_{\mathbf{i}}\leq u_{\mathbf{n,i}})=O(P(X_{\mathbf{j}}\leq u_{\mathbf{n,j}}))$ and $P(X_{\mathbf{i}}\leq v_{\mathbf{n,i}})=O(P(X_{\mathbf{j}}\leq v_{\mathbf{n,j}}))$ for all $\mathbf{1\leq i\neq j\leq n}$ as $\mathbf{n}\rightarrow\infty$. This makes the random field $\mathbf{X}$ look like a stationary one. However, this condition is necessary to deal with the random missing data (see Remark 2.1 for details).
Let $\mathbf{{R_\mathbf{n}}}=\{1,\cdots,n_1\}\times\{1,\cdots,n_2\}$ and $\mathbf{1}=(1,1)$,  and subdivide $\mathbf{R_n}$ into $k_{n_{1}}k_{n_{2}}$ disjoint rectangle subsets $\mathbf{K_{s}}=\mathbf{K}_{(s_{1},s_{2})}, s_{1}=1,2,\ldots, k_{n_{1}}, s_{2}=1,2,\ldots, k_{n_{2}}$, such that
\begin{equation}
\label{eqTan1}
\sum_{i\in \mathbf{K_{s}}}P\left(X_\mathbf{i}>w_{\mathbf{n,i}}\right)=\frac{1}{k_{n_1}k_{n_2}}\sum_{\mathbf{i}\in \mathbf{R_{n}}}P\left(X_{\mathbf{i}}>w_{\mathbf{n,i}}\right)+o(1),
\end{equation}
as $\mathbf{n}\rightarrow\infty$, where $w_{\mathbf{n,i}}$ equals $v_{\mathbf{n,i}}$ or $u_{\mathbf{n,i}}$ for all $\mathbf{i}\leq \mathbf{n}$.
Note that average partitions of $\mathbf{{R_\mathbf{n}}}$ satisfies (\ref{eqTan1}) under the above assumptions.

\textbf{Definition 1.1.} Let $\mathcal{F}$ be a family of rectangle subsets $\mathbf{K_{s}}$. The nonstationary random field $\mathbf{X}$ on $\mathbb{Z}_{+}^2$ satisfies the condition $\mathbf{D}(u_{\mathbf{n,i}}, v_{\mathbf{n,i}})$
over $\mathcal{F}$ if there exist sequences of integer valued constants $\{k_{n_{i}}\}_{n_{i}\geq 1}$,$\{m_{n_{i}}\}_{n_{i}\geq 1},i=1,2$ with \\
\begin{equation}
\label{eqTan2}
\left(k_{n_{1}},k_{n_{2}}\right)\rightarrow \infty,\ \
\left(\frac{k_{n_1}m_{n_1}}{n_1},\frac{k_{n_2}m_{n_2}}{n_2}\right)\rightarrow
\mathbf{0}
\end{equation}
such that $k_{n_1}k_{n_2}\alpha_{\mathbf{n},m_{n_1},m_{n_2}}\rightarrow 0$, as $\mathbf{n}=(n_{1},n_{2})\rightarrow \mathbf{\infty}$,
 where $\alpha_{\mathbf{n},m_{n_1},m_{n_2}}$ is defined as
\begin{equation}
\begin{split}
\label{eqTan3}
\alpha_{\mathbf{n},m_{n_1},m_{n_2}}=\sup_{\substack{(\mathbf{I}_1\bigcup \mathbf{J}_1, \mathbf{I}_2 \bigcup \mathbf{J}_2)\\ \in \mathcal{S}(m_{n_1},m_{n_2})}}\bigg|P\left(\bigcap_{\mathbf{i}\in \mathbf{I}_1\bigcup \mathbf{J}_1}\{X_{\mathbf{i}}\leq u_{\mathbf{n,i}}\}\cap \bigcap_{\mathbf{j}\in \mathbf{I}_2\bigcup \mathbf{J}_2}\{X_{\mathbf{j}} \leq v_{\mathbf{n,j}}\}\right)-\\P\left(\bigcap_{\mathbf{i}\in \mathbf{I}_1}\{X_\mathbf{i}\leq u_{\mathbf{n,i}}\}\cap  \bigcap_{\mathbf{j}\in \mathbf{I}_2}\{X_{\mathbf{j}}\leq v_{\mathbf{n,j}}\}\right)P\left(\bigcap_{\mathbf{i}\in \mathbf{J}_1}\{X_\mathbf{i}\leq u_{\mathbf{n,i}}\}\cap\bigcap_{\mathbf{j}\in \mathbf{J}_2}\{X_j\leq v_{\mathbf{n,j}}\}\right)\bigg|
\end{split}
\end{equation}
 with $\mathcal{S}(m_{n_1},m_{n_2})=\{(\mathbf{I},\mathbf{J})\in \mathcal{F}^2,d(\pi_1(\mathbf{I}),\pi_1(\mathbf{J}))\geq m_{n_1}\bigvee d(\pi_2(\mathbf{I}),\pi_2(\mathbf{J}))\geq m_{n_2}\}$, $\mathbf{I}_1\bigcap \mathbf{I}_2=\varnothing, \mathbf{J}_1\bigcap \mathbf{J}_2=\varnothing$. Here $\Pi_{i}, i=1,2$, denote the cartesian projections on $x$-axis and $y$-axis, and $d(A,B)$ denotes the distance between two sets $A$ and $B$.

The next local dependent condition is taken from \cite{Pereira_Tan2017}.

\textbf{Definition 1.2.}
The condition $\mathbf{D}'(v_{\mathbf{n,i}})$ holds for the nonstationary random field $\mathbf{X}$,
if for each $\mathbf{I}\in\mathcal{F}$ , we have
\begin{equation}
\label{eqTan4}
k_{n_1}k_{n_2}\sum_{\mathbf{i\neq j}\in \mathbf{I}}P\left(X_{\mathbf{i}}>v_{\mathbf{n,i}},X_{\mathbf{j}}>v_{\mathbf{n,j}}\right)\rightarrow 0,\ \ as\ \
\mathbf{n}\rightarrow \mathbf{\infty}.
\end{equation}
The local dependence condition $\mathbf{D}'(v_{\mathbf{n, i}})$ is an anti-cluster condition, which bounds the probability of more than one exceedance above the levels $v_{\mathbf{n,i}}$ in a rectangle with a few indexes.

\section{Main results}

Assume that some of random variables in the random filed $\mathbf{X}$ can be observed and the sequence of random variables $\varepsilon=\{\varepsilon_{\mathbf{i}},\mathbf{i}\geq \mathbf{1}\}$ indicate which variables in the random filed $\mathbf{X}$ are observed.
Let $S_{\mathbf{n}}=\sum_{\mathbf{i}\leq \mathbf{n}}\varepsilon_{\mathbf{i}}$ and suppose that
\begin{equation}
\label{eqT11.5}
\frac{S_{\mathbf{n}}}{n_{1}n_{2}}\stackrel{P}{\longrightarrow} \lambda,\ \ \mbox{as}\ \ \mathbf{n}\rightarrow\infty,
\end{equation}
where $\lambda$ is a random variable satisfying $0\leq \lambda\leq 1$ a.s. For any random variable sequence $\{\xi_{\mathbf{i}}, \mathbf{n\geq 1}\}$, we  define
\begin{equation}
\label{eqT11.6}
\xi_{\mathbf{i}}(\varepsilon)=(1-\varepsilon_{\mathbf{i}})\gamma(\xi_{\mathbf{i}})+\varepsilon_{\mathbf{i}}\xi_{\mathbf{i}},\ \ \mathbf{i\geq 1},
\end{equation}
where
$\gamma(\xi_{\mathbf{i}})=\inf\{x\in \mathbb{R}: P(\xi_{\mathbf{i}}\leq x)>0\}$.

\textbf{Theorem 2.1}. {\sl Suppose that the nonstationary random field $\mathbf{X}=\{X_{\mathbf{n}},\mathbf{n}\geq \mathbf{1}\}$ satisfies $\mathbf{D}(u_{\mathbf{n,i}},v_{\mathbf{n,i}})$ and
$\mathbf{D}'(v_{\mathbf{n,i}})$ over $\mathcal{F}$ and
$\sup_{\mathbf{n}\geq\mathbf{1}}\{n_{1}n_{2}P(X_{\mathbf{i}}\geq v_{\mathbf{n},\mathbf{i}}), \mathbf{i}\leq \mathbf{n}\}$ is bounded.
Assume that $\mathbf{\varepsilon}=\{\varepsilon_{\mathbf{n}},\mathbf{n}\geq \mathbf{1}\}$ is a sequence of indicators that is independent of $\mathbf{X}=\{X_{\mathbf{n}},\mathbf{n}\geq \mathbf{1}\}$ and such that (\ref{eqT11.5}) holds.
If $\sum_{\mathbf{i}\in \mathbf{R_{n}}}P\left(X_\mathbf{i}>u_{\mathbf{n,i}}\right)\rightarrow\tau>0$ and $\sum_{\mathbf{i}\in \mathbf{R_{n}}}P\left(X_\mathbf{i}>v_{\mathbf{n,i}}\right)\rightarrow\kappa>0$ hold,
then
\begin{equation}
\label{eqT11.7}
\lim_{\mathbf{n}\rightarrow \mathbf{\infty}}P\left(\bigcap_{\mathbf{i}\in \mathbf{R_{n}}}\{X_{\mathbf{i}}(\varepsilon)\leq v_{\mathbf{n,i}}\},\bigcap_{\mathbf{i}\in \mathbf{R_{n}}}\{X_{\mathbf{i}}\leq u_{\mathbf{n,i}}\}\right)=E[e^{-\lambda\kappa}e^{-(1-\lambda)\tau}].
\end{equation}
}

\textbf{Remark 2.1}.  Let us divide $\mathbf{R_{n}}$ onto two parts $\mathbf{R}_{\mathbf{n},1}$ and $\mathbf{R}_{\mathbf{n},2}$ and let $P(X_{\mathbf{i}}\leq u_{\mathbf{n,i}})=o(P(X_{\mathbf{j}}\leq u_{\mathbf{n,j}}))$  and $P(X_{\mathbf{i}}\leq v_{\mathbf{n,i}})=o(P(X_{\mathbf{j}}\leq v_{\mathbf{n,j}}))$ for $\mathbf{i}\in\mathbf{R}_{\mathbf{n},1}$ and $\mathbf{j}\in\mathbf{R}_{\mathbf{n},2}$, respectively. Thus, $\sum_{\mathbf{i}\in \mathbf{R}_{\mathbf{n},2}}P\left(X_\mathbf{i}>u_{\mathbf{n,i}}\right)\rightarrow\tau>0$ and $\sum_{\mathbf{i}\in \mathbf{R}_{\mathbf{n},2}}P\left(X_\mathbf{i}>v_{\mathbf{n,i}}\right)\rightarrow\kappa>0$. Since we do not know the limit of $\frac{S_{\mathbf{R}_{\mathbf{n},2}}}{n_{1}n_{2}}$, the limit of the probability in (\ref{eqT11.7}) does not exist,
where $S_{\mathbf{R}_{\mathbf{n},2}}=\sum_{\mathbf{i}\in \mathbf{R}_{\mathbf{n},2}}\varepsilon_{\mathbf{i}}$. This shows that
the condition that $P(X_{\mathbf{i}}\leq u_{\mathbf{n,i}})=O(P(X_{\mathbf{j}}\leq u_{\mathbf{n,j}}))$ and $P(X_{\mathbf{i}}\leq v_{\mathbf{n,i}})=O(P(X_{\mathbf{j}}\leq v_{\mathbf{n,j}}))$ for all $\mathbf{1\leq i\neq j\leq n}$ as $\mathbf{n}\rightarrow\infty$ is necessary.

Next, we extend the above weakly convergence results to almost sure version. As usual, $a_{\mathbf{n}}\ll b_{\mathbf{n}}$ means $a_{\mathbf{n}}=O(b_{\mathbf{n}})$.
In order to formulate the results, we need another condition $\mathbf{D}^{*}(u_{\mathbf{k,j}},v_{\mathbf{k,i}}, u_{\mathbf{n,j}},v_{\mathbf{n,i}})$ as follows.

\textbf{Definition 2.1.} The random field $\mathbf{X}$ on $\mathbb{Z}_{+}^2$ satisfies the condition $\mathbf{D}^{*}(u_{\mathbf{k,j}},v_{\mathbf{k,i}}, u_{\mathbf{n,j}},v_{\mathbf{n,i}})$ if there exist sequences of integer valued constants $\{m_{n_{i}}\}_{n_{i}\geq 1}, i=1, 2$,  and
\begin{eqnarray*}
\left(m_{n_{1}},m_{n_{2}}\right)\rightarrow \infty,\ \ \left(\frac{m_{n_1}}{n_1},\frac{m_{n_2}}{n_2}\right)\rightarrow
\mathbf{0}
\end{eqnarray*}
such that for some $\epsilon>0$
\begin{equation}
\sup_{k_{1}k_{2}< n_{1}n_{2}}\alpha_{\mathbf{n},\mathbf{k},m_{n_1},m_{n_2}}^{*}\ll \left(\log\log n_1\log\log n_2\right)^{-(1+\epsilon)}
\end{equation}
as $\mathbf{n}=(n_{1},n_{2})\rightarrow \mathbf{\infty}$, where $\alpha_{\mathbf{n},\mathbf{k},m_{n_1},m_{n_2}}^{*}$ is defined as following. For any  $\mathbf{k\neq n}$ such that $k_{1}k_{2}< n_{1}n_{2}$ define
\begin{eqnarray*}
&&\alpha_{\mathbf{n},\mathbf{k},m_{n_1},m_{n_2}}^{*}\\
&&= \sup_{\mathbf{I}_{1}\subseteq \mathbf{R_{k}}, \mathbf{I}_{2}\subseteq\mathbf{R_n}\backslash \mathbf{M_{kn}}}\bigg|P(\bigcap_{\mathbf{i}\in \mathbf{I}_1}\{X_{\mathbf{i}}\leq v_{\mathbf{k},\mathbf{i}}\}\cap \bigcap_{\mathbf{j} \in \mathbf{J}_1}\{X_{\mathbf{j}}\leq u_{\mathbf{k},\mathbf{j}}\}\cap\bigcap_{\mathbf{i}\in \mathbf{I}_2}\{X_{\mathbf{i}}\leq v_{\mathbf{n},\mathbf{i}}\}\cap\bigcap_{\mathbf{j} \in \mathbf{J}_2}\{X_{\mathbf{j}}\leq u_{\mathbf{n},\mathbf{j}}\})\\
&&\ \ \ \ \ \ \ \ \ \ \ \ \ \ \ \ \ \ \ \ \ \ \ \ -P(\bigcap_{\mathbf{i}\in \mathbf{I}_1}\{X_{\mathbf{i}}\leq v_{\mathbf{k},\mathbf{i}}\}\cap \bigcap_{\mathbf{j} \in \mathbf{J}_1}\{X_{\mathbf{j}}\leq u_{\mathbf{k},\mathbf{j}}\})P(\bigcap_{\mathbf{i}\in \mathbf{I}_2}\{X_{\mathbf{i}}\leq v_{\mathbf{n},\mathbf{i}}\}\cap\bigcap_{\mathbf{j} \in \mathbf{J}_2}\{X_{\mathbf{j}}\leq u_{\mathbf{n},\mathbf{j}}\})\bigg|
\end{eqnarray*}
where $\mathbf{I}_1\subseteq \mathbf{J}_1=\mathbf{R_k}, \mathbf{I}_2\subseteq \mathbf{J}_2=\mathbf{R_n}\backslash \mathbf{M_{kn}}$ with $\mathbf{M_{kn}}=\{(j_1,j_2):(j_1,j_2)\in \mathbb{N}^2,0 \leq j_i \leq \#(\pi_i(\mathbf{M^{*}_{kn}}))+m_{n_i},i=1,2\}$ and $ \mathbf{M}_{\mathbf{kn}}^{*}
=\mathbf{R_k}\bigcap \mathbf{R_n}$. Here $\#(A)$ denotes cardinality of the set $A$.

\textbf{Theorem 2.2}. {\sl Suppose that the nonstationary random field $\mathbf{X}=\{X_{\mathbf{n}},\mathbf{n}\geq \mathbf{1}\}$ satisfies $\mathbf{D}(u_{\mathbf{n,i}},v_{\mathbf{n,i}})$, $\mathbf{D}^{*}(u_{\mathbf{k,j}},v_{\mathbf{k,i}}, u_{\mathbf{n,j}},v_{\mathbf{n,i}})$ and
$\mathbf{D}'(v_{\mathbf{n,i}})$ over $\mathcal{F}$ and
$\sup_{\mathbf{n}\geq\mathbf{1}}\{n_{1}n_{2}P(X_{\mathbf{i}}\geq v_{\mathbf{n},\mathbf{i}}), \mathbf{i}\leq \mathbf{n}\}$ is bounded.
Assume that $\mathbf{\varepsilon}=\{\varepsilon_{\mathbf{n}},\mathbf{n}\geq \mathbf{1}\}$ is a sequence of independent indicators that is independent of $\mathbf{X}=\{X_{\mathbf{n}},\mathbf{n}\geq \mathbf{1}\}$ and such that (\ref{eqT11.5}) holds.
If $\sum_{\mathbf{i}\in \mathbf{R_{n}}}P\left(X_\mathbf{i}>u_{\mathbf{n,i}}\right)\rightarrow\tau>0,\sum_{\mathbf{i}\in \mathbf{R_{n}}}P\left(X_\mathbf{i}>v_{\mathbf{n,i}}\right)\rightarrow\kappa>0$, and $n_{1}=O(n_{2})$ hold,
then
\begin{equation}
\label{eqT21.2}
\lim_{\mathbf{n}\rightarrow \infty}\frac{1}{\log n_1 \log n_2}\sum_{\mathbf{k}\in \mathbf{R_n}}\frac{1}{k_1k_2}\mathbf{\mathbbm{1}}_{\left\{\bigcap_{\mathbf{i}\in \mathbf{R_k}} \{X_{\mathbf{i}}(\varepsilon)\leq v_{\mathbf{k,i}}, X_{\mathbf{i}}\leq u_{\mathbf{k,i}}\}\right\}}= E[e^{-\lambda\kappa}e^{-(1-\lambda)\tau}],\ \  a.s.
\end{equation}
}

\textbf{Remark 2.2}. Theorem 2.1 of \cite{Pereira_Tan2017} derived the almost sure limit theorem for the maxima of nonstationary random fields, but the condition $\mathbf{D}^{*}(u_{\mathbf{n,i}})$ in their paper can not imply their main result. More precisely, for the case $k_{1}\leq l_{1}, k_{2}\leq l_{2}$,
the term $I_{2}$ in the proof of Lemma 4.1 in their paper can not be bounded by $\alpha_{\mathbf{l},m_{l_{1}},m_{l_{2}}}$, since the sets $\mathbf{R_{k}}$ and $\mathbf{R_{l}-M_{kl}}$ do not satisfy the conditions of $\mathbf{D}^{*}(u_{\mathbf{n,i}})$. We need a technical condition $\mathbf{D}^{*}(u_{\mathbf{k,j}},v_{\mathbf{k,i}}, u_{\mathbf{n,j}},v_{\mathbf{n,i}})$ to solve this problem.  However this condition is not very strict but involve more levels and different index sets. We will give several examples to  illustrate it.

We end with this section by two examples which satisfy the conditions of Theorems 2.1 and 2.2. For examples such as for Gaussian case and its function, it will be given as applications in Section 3.

\textbf{Example 2.1.} {\sl  Suppose that $\mathbf{X}=\{X_{\mathbf{n}},\mathbf{n}\geq \mathbf{1}\}$ is a sequence of independent
random fields, then conditions $\mathbf{D}(u_{\mathbf{n,i}},v_{\mathbf{n,i}})$, $\mathbf{D}^{*}(u_{\mathbf{k,j}},v_{\mathbf{k,i}}, u_{\mathbf{n,j}},v_{\mathbf{n,i}})$ and $\mathbf{D}'(v_{\mathbf{n,i}})$ hold.
}

\textbf{Example 2.2.} {\sl  Suppose that $\mathbf{X}=\{X_{\mathbf{n}},\mathbf{n}\geq \mathbf{1}\}$ is a sequence of $m-$denpendent or strong mixing random fields, then condition $\mathbf{D}(u_{\mathbf{n,i}},v_{\mathbf{n,i}})$ holds.
}

\section{Applications}

Applying our main results, we can derive the weak and almost sure convergence results of the maxima for complete and incomplete samples for nonstationary Gaussian random fields and their functions.

Suppose in this section that $\mathbf{Y}=\{Y_{\mathbf{n}},\mathbf{n}\geq \mathbf{1}\}$ is a sequence of standard (with 0 mean and 1 variance) Gaussian random fields with covariances $r_{\mathbf{i,j}}=Cov(Y_{\mathbf{i}},Y_{\mathbf{j}})$.
Assume that $r_{\mathbf{i,j}}$ satisfies $|r_{\mathbf{i,j}}|<\rho_{|\mathbf{i-j}|}$ when $\mathbf{i\neq j}$, for some sequence $\rho_{\mathbf{n}}<1$ for $\mathbf{n}\neq \mathbf{0}$ such that
\begin{eqnarray}
\label{21}
\rho_{(n_{1},0)}\log n_{1}\ \mbox{and}\ \rho_{(0,n_{2})}\log n_{2}\ \mbox{are\ \ bounded}\ \mbox{and}\ \ \lim_{\mathbf{n}\rightarrow\infty}\rho_{\mathbf{n}}\log(n_{1}n_{2})=0
\end{eqnarray}
or
\begin{eqnarray}
\label{t21}
&&\rho_{(n_{1},0)}\log n_{1}\ll (\log\log n_{1})^{-(1+\epsilon)},\ \rho_{(0,n_{2})}\log n_{2}\ll (\log\log n_{2})^{-(1+\epsilon)}\ \ \mbox{and}\nonumber\\
&& \ \rho_{\mathbf{n}}\log(n_{1}n_{2})\ll (\log\log n_{1}\log\log n_{2})^{-(1+\epsilon)}.
\end{eqnarray}
Let $\{u_{\mathbf{n,i}}, \mathbf{i\leq n}\}_{\mathbf{n \geq 1}}$ and $\{v_{\mathbf{n,i}},\mathbf{i\leq n}\}_{\mathbf{n \geq 1}}$ be two sequences of real numbers satisfying that $v_{\mathbf{n,i}}\leq u_{\mathbf{n,i}}$ for all $\mathbf{i\leq n}$ and $\Phi(u_{\mathbf{n,i}})=O(\Phi(u_{\mathbf{n,j}}))$, $\Phi(v_{\mathbf{n,i}})=O(\Phi(v_{\mathbf{n,j}}))$ for all $\mathbf{1\leq i\neq j\leq n}$ and sufficiently large $\mathbf{n}$.

\textbf{Theorem 3.1.} {\sl  Suppose that $\mathbf{Y}=\{Y_{\mathbf{n}},\mathbf{n}\geq \mathbf{1}\}$ is a sequence of standard Gaussian random fields with covariances satisfying (\ref{21}). Assume that $\mathbf{\varepsilon}=\{\varepsilon_{\mathbf{n}},\mathbf{n}\geq \mathbf{1}\}$ is a sequence of indicators that is independent of $\mathbf{Y}$ and such that (\ref{eqT11.5}) holds.
Let the constants $\{u_{\mathbf{n,i}}, \mathbf{i\leq n}\}_{\mathbf{n \geq 1}}$ and
$\{v_{\mathbf{n,i}},\mathbf{i\leq n}\}_{\mathbf{n \geq 1}}$ be such that $\sup_{\mathbf{n}\geq\mathbf{1}}\{n_{1}n_{2}(1-\Phi(v_{\mathbf{n},\mathbf{i}}), \mathbf{i}\leq \mathbf{n})\}$ is bounded, $\sum_{\mathbf{i}\in \mathbf{R_{n}}}[1-\Phi(u_{\mathbf{n,i}})]\rightarrow\tau>0$ and $\sum_{\mathbf{i}\in \mathbf{R_{n}}}[1-\Phi(v_{\mathbf{n,i}})]\rightarrow\kappa>0$.
Then, we have
$$\lim_{\mathbf{n}\rightarrow \mathbf{\infty}}P\left(\bigcap_{\mathbf{i}\in \mathbf{R_{n}}}\{Y_{\mathbf{i}}(\varepsilon)\leq v_{\mathbf{n,i}}\},\bigcap_{\mathbf{i}\in \mathbf{R_{n}}}\{Y_{\mathbf{i}}\leq u_{\mathbf{n,i}}\}\right)=E[e^{-\lambda\kappa}e^{-(1-\lambda)\tau}].$$
Furthermore, if $\mathbf{\varepsilon}=\{\varepsilon_{\mathbf{n}},\mathbf{n}\geq \mathbf{1}\}$ is a sequence of independent indicators, $n_{1}=O(n_{2})$ and (\ref{t21}) holds, we have
$$
\lim_{\mathbf{n}\rightarrow \infty}\frac{1}{\log n_1 \log n_2}\sum_{\mathbf{k}\in \mathbf{R_n}}\frac{1}{k_1k_2}\mathbf{\mathbbm{1}}_{\left\{\bigcap_{\mathbf{i}\in \mathbf{R_k}} \{Y_{\mathbf{i}}(\varepsilon)\leq v_{\mathbf{k,i}}, Y_{\mathbf{i}}\leq u_{\mathbf{k,i}}\}\right\}}= E[e^{-\lambda\kappa}e^{-(1-\lambda)\tau}],\ \  a.s.
$$
}

\textbf{Remark 3.1.} 1) Under the same conditions as Theorem 3.1, \cite{Tan_Wang2014} derived the almost sure limit theorem for maxima of complete samples of nonstationary Gaussian random fields. Theorem 3.1 extended their results to the maxima for complete and incomplete samples.\\
2) For one-dimensional case, i.e., for Gaussian sequences, the condition that $\mathbf{\varepsilon}=\{\varepsilon_{n},n\geq 1\}$ is a sequence of independent indicators can be weakened as that $\mathbf{\varepsilon}=\{\varepsilon_{n},n\geq 1\}$ is strong mixing indicators with mixing coefficient $\alpha(n)\ll (\log\log n)^{-(1+\epsilon)}$ for some $\epsilon>0$.\\
3) The one-dimensional case of Theorem 3.1 is an extension of the main results of \cite{Tong_Peng2011} which derived the almost limit theorem of
the maxima for complete and incomplete samples for stationary Gaussian case when $\lambda$ is a constant.

Let $a_{\mathbf{n}}=\sqrt{2\log(n_{1}n_{2})}$ and $b_{\mathbf{n}}=a_{\mathbf{n}}-\frac{\log\log(n_{1}n_{2})+\log (4\pi)}{2a_{\mathbf{n}}}$.

\textbf{Corollary 3.1}. {\sl Let $\mathbf{Z}=\{Y_{\mathbf{n}}+m_{\mathbf{n}},\mathbf{n}\geq \mathbf{1}\}$, where $\{Y_{\mathbf{n}},\mathbf{n}\geq \mathbf{1}\}$
is defined as above and $\{m_{\mathbf{n}},\mathbf{n}\geq \mathbf{1}\}$
satisfies
\begin{equation}
\label{Tan2.7}
\beta_{\mathbf{n}}=\max_{\mathbf{k}\in \mathbf{R_{n}}}|m_{\mathbf{k}}|=o(\sqrt{n_{1}n_{2}}),
\end{equation}
and let $m_{\mathbf{n}}^{*}$ be such that
\begin{equation}
\label{Tan2.8}
|m_{\mathbf{n}}^{*}|\leq\beta_{\mathbf{n}}
\end{equation}
and
\begin{equation}
\label{Tan2.9}
\frac{1}{n_{1}n_{2}}\sum_{\mathbf{i}\in \mathbf{R_{n}}}\exp\bigg(a_{\mathbf{n}}^{*}(m_{\mathbf{i}}-m_{\mathbf{n}}^{*})-\frac{1}{2}(m_{\mathbf{i}}-m_{\mathbf{n}}^{*})^{2}\bigg)\rightarrow 1
\end{equation}
as $\mathbf{n}\rightarrow\infty$, where $a_{\mathbf{n}}^{*}=a_{\mathbf{n}}-\log\log(n_{1}n_{2})/2a_{\mathbf{n}}$.
Assume that $\mathbf{\varepsilon}=\{\varepsilon_{\mathbf{n}},\mathbf{n}\geq \mathbf{1}\}$ is a sequence of indicators that is independent of $\mathbf{Y}$ and such that (\ref{eqT11.5}) holds.
If (\ref{21}) holds, then for any $x\leq y\in \mathbb{R}$
\begin{eqnarray}
\label{Tan2.10}
&&\lim_{\mathbf{n}\rightarrow\infty}P\bigg(a_{\mathbf{n}}\big(M_{\mathbf{n}}(Z(\varepsilon))-b_{\mathbf{n}}-m_{\mathbf{n}}^{*}\big)\leq x, a_{\mathbf{n}}\big(M_{\mathbf{n}}(Z)-b_{\mathbf{n}}-m_{\mathbf{n}}^{*}\big)\leq y\bigg)\nonumber\\
&&=E[\exp(-\lambda e^{-x})\exp(-(1-\lambda) e^{-y})];
\end{eqnarray}
Furthermore, if $\mathbf{\varepsilon}=\{\varepsilon_{\mathbf{n}},\mathbf{n}\geq \mathbf{1}\}$ is a sequence of independent indicators, $n_{1}=O(n_{2})$, (\ref{t21}) holds and
\begin{equation}
\label{eq2.11}
a_{\mathbf{n}}\big(\max_{\mathbf{i}\in \mathbf{R_{n}}}m_{\mathbf{i}}-m_{\mathbf{n}}^{*}\big)\ \ \mbox{is\ \ bounded},
\end{equation}
 then for any $x\leq y\in \mathbb{R}$
\begin{eqnarray}
\label{Tan2.10}
&&\lim_{\mathbf{n}\rightarrow \infty}\frac{1}{\log n_1 \log n_2}\sum_{\mathbf{k}\in \mathbf{R_n}}\frac{1}{k_1k_2}\mathbf{\mathbbm{1}}_{\left\{a_{\mathbf{k}}\big(M_{\mathbf{k}}(Z(\varepsilon))-b_{\mathbf{k}}-m_{\mathbf{k}}^{*}\big)\leq x, a_{\mathbf{k}}\big(M_{\mathbf{k}}(Z)-b_{\mathbf{k}}-m_{\mathbf{k}}^{*}\big)\leq y\right\}}\nonumber\\
&&=E[\exp(-\lambda e^{-x})\exp(-(1-\lambda) e^{-y})],\ \  a.s.,
\end{eqnarray}
where $M_{\mathbf{n}}(Z(\varepsilon))=\max_{\mathbf{i}\in \mathbf{R_{n}}}Z_{\mathbf{i}}(\varepsilon)$ and $M_{\mathbf{n}}(Z)=\max_{\mathbf{i}\in \mathbf{R_{n}}}Z_{\mathbf{i}}$.
}

Next, we deal with two types of Gaussian functions. Let $\{Y_{\mathbf{n}j},\mathbf{n}\geq \mathbf{1}\}, j=1,2,\cdots,d$, $d\geq 1$ be the independent copies of $\{Y_{\mathbf{n}},\mathbf{n}\geq \mathbf{1}\}$.
Define
$$\chi_{\mathbf{n}}=(\sum_{j=1}^{d}Y_{\mathbf{n}j}^{2})^{1/2},\ \ \mathbf{n}\geq \mathbf{1}$$
and  for
$r\in\{1,2,\ldots,d\}$
$$O_{\mathbf{n}}^{(d)}:=\min_{j=1}^{d}Y_{\mathbf{n}j}\leq \cdots\leq O_{\mathbf{n}}^{(r)}\leq \cdots\leq O_{\mathbf{n}j}^{(1)}:=\max_{j=1}^{d}Y_{\mathbf{n}j},\ \ \mathbf{n}\geq \mathbf{1}.$$
It worth pointing out that $\{\chi_{\mathbf{n}},\mathbf{n}\geq \mathbf{1}\}$ is a sequence of $\chi$ random field and $\{O_{\mathbf{n}}^{(r)},\mathbf{n}\geq \mathbf{1}\}$ is a
sequence of Gaussian order statistic random field. We refer to \cite{Tan_Hashorva2013a, Tan_Hashorva2013b}, \cite{Ling_Tan2016} and \cite{Shao_Tan2022} for recent work on extremes for $\chi$ variables and \cite{Debicki_Hashorva_Ji_Ling2015, Debicki_Hashorva_Ji_Ling2017} and \cite{Tan2018} for recent work on extremes for Gaussian order statistic variables.

\textbf{Theorem 3.2.} {\sl  Suppose that $\{Y_{\mathbf{n}},\mathbf{n}\geq \mathbf{1}\}$ is a sequence of standard Gaussian random fields with covariances satisfying (\ref{21}). Assume that $\mathbf{\varepsilon}=\{\varepsilon_{\mathbf{n}},\mathbf{n}\geq \mathbf{1}\}$ is a sequence of indicators that is independent of $\mathbf{Y}$ and such that (\ref{eqT11.5}) holds.
Let the constants $\{u_{\mathbf{n,i}}, \mathbf{i\leq n}\}_{\mathbf{n \geq 1}}$ and
$\{v_{\mathbf{n,i}},\mathbf{i\leq n}\}_{\mathbf{n \geq 1}}$ be such that $\sup_{\mathbf{n}\geq\mathbf{1}}\{n_{1}n_{2}P(\chi_\mathbf{i}>v_{\mathbf{n},\mathbf{i}}), \mathbf{i}\leq \mathbf{n})\}$ is bounded, $\sum_{\mathbf{i}\in \mathbf{R_{n}}}P\left(\chi_\mathbf{i}>u_{\mathbf{n,i}}\right)\rightarrow\tau>0$ and $\sum_{\mathbf{i}\in \mathbf{R_{n}}}P\left(\chi_\mathbf{i}>v_{\mathbf{n,i}}\right)\rightarrow\kappa>0$ hold. Then,
$$\lim_{\mathbf{n}\rightarrow \mathbf{\infty}}P\left(\bigcap_{\mathbf{i}\in \mathbf{R_{n}}}\{\chi_{\mathbf{i}}(\varepsilon)\leq v_{\mathbf{n,i}}\},\bigcap_{\mathbf{i}\in \mathbf{R_{n}}}\{\chi_{\mathbf{i}}\leq u_{\mathbf{n,i}}\}\right)=E[e^{-\lambda\kappa}e^{-(1-\lambda)\tau}].$$
Furthermore,  if $\mathbf{\varepsilon}=\{\varepsilon_{\mathbf{n}},\mathbf{n}\geq \mathbf{1}\}$ is a sequence of independent indicators, $n_{1}=O(n_{2})$ and (\ref{t21}) holds, we have
$$
\lim_{\mathbf{n}\rightarrow \infty}\frac{1}{\log n_1 \log n_2}\sum_{\mathbf{k}\in \mathbf{R_n}}\frac{1}{k_1k_2}\mathbf{\mathbbm{1}}_{\left\{\bigcap_{\mathbf{i}\in \mathbf{R_k}} \{\chi_{\mathbf{i}}(\varepsilon)\leq v_{\mathbf{k,i}}, \chi_{\mathbf{i}}\leq u_{\mathbf{k,i}}\}\right\}}= E[e^{-\lambda\kappa}e^{-(1-\lambda)\tau}],\ \  a.s.
$$
}

\textbf{Theorem 3.3.} {\sl  Suppose that $\{Y_{\mathbf{n}},\mathbf{n}\geq \mathbf{1}\}$ is a sequence of standard Gaussian random fields with covariances satisfying (\ref{21}). Assume that $\mathbf{\varepsilon}=\{\varepsilon_{\mathbf{n}},\mathbf{n}\geq \mathbf{1}\}$ is a sequence of indicators that is independent of $\mathbf{Y}$ and such that (\ref{eqT11.5}) holds.
Let the constants $\{u_{\mathbf{n,i}}, \mathbf{i\leq n}\}_{\mathbf{n \geq 1}}$ and
$\{v_{\mathbf{n,i}},\mathbf{i\leq n}\}_{\mathbf{n \geq 1}}$ be such that $\sup_{\mathbf{n}\geq\mathbf{1}}\{n_{1}n_{2}P(O_{\mathbf{i}}^{(r)}>v_{\mathbf{n},\mathbf{i}}), \mathbf{i}\leq \mathbf{n})\}$ is bounded, $\sum_{\mathbf{i}\in \mathbf{R_{n}}}P\left(O_{\mathbf{i}}^{(r)}>u_{\mathbf{n,i}}\right)\rightarrow\tau>0$ and $\sum_{\mathbf{i}\in \mathbf{R_{n}}}P\left(O_{\mathbf{i}}^{(r)}>v_{\mathbf{n,i}}\right)\rightarrow\kappa>0$ hold. Then,
$$\lim_{\mathbf{n}\rightarrow \mathbf{\infty}}P\left(\bigcap_{\mathbf{i}\in \mathbf{R_{n}}}\{O_{\mathbf{i}}^{(r)}(\varepsilon)\leq v_{\mathbf{n,i}}\},\bigcap_{\mathbf{i}\in \mathbf{R_{n}}}\{O_{\mathbf{i}}^{(r)}\leq u_{\mathbf{n,i}}\}\right)=E[e^{-\lambda\kappa}e^{-(1-\lambda)\tau}].$$
Furthermore,  if $\mathbf{\varepsilon}=\{\varepsilon_{\mathbf{n}},\mathbf{n}\geq \mathbf{1}\}$ is a sequence of independent indicators, $n_{1}=O(n_{2})$ and (\ref{t21}) holds, we have
$$
\lim_{\mathbf{n}\rightarrow \infty}\frac{1}{\log n_1 \log n_2}\sum_{\mathbf{k}\in \mathbf{R_n}}\frac{1}{k_1k_2}\mathbf{\mathbbm{1}}_{\left\{\bigcap_{\mathbf{i}\in \mathbf{R_k}} \{O_{\mathbf{i}}^{(r)}(\varepsilon)\leq v_{\mathbf{k,i}}, O_{\mathbf{i}}^{(r)}\leq u_{\mathbf{k,i}}\}\right\}}= E[e^{-\lambda\kappa}e^{-(1-\lambda)\tau}],\ \  a.s.
$$
}

\section{Auxiliary results and proofs}

In this section, we first state and prove several lemmas which will be used in the proofs of our main results and then we give the proofs  of the main results.
For any $\mathbf{I}\subseteq \mathbf{R_{n}}$, let $\mathcal{B}_\mathbf{k}(\mathbf{I})=\bigcap_{\mathbf{i}\in {\mathbf{I}}}\{X_{\mathbf{i}}\leq u_{\mathbf{k,i}}, X_{\mathbf{i}}(\varepsilon)\leq v_{\mathbf{k,i}}\}$, and  $\overline{\mathcal{B}}_\mathbf{k}(\mathbf{I})=\bigcup_{\mathbf{i}\in {\mathbf{I}}}\{\{X_{\mathbf{i}}> u_{\mathbf{k,i}}\}\bigcup \{X_{\mathbf{i}}(\varepsilon)> v_{\mathbf{k,i}}\}\}$. Let $m_{l_i}=\log l_i$, $i=1,2$.
For random variable $\lambda$
 such that $0\leq \lambda \leq 1$ a.s., we put
 \begin{eqnarray*}
 &&\ \ B_{r,\mathbf{l}}=\left\{
   \begin{array}{ll}
   \omega:\lambda(\omega)\in
      \left\{
		\begin{array}{ll}
			\left[0,~~\frac{1}{2^{l_1l_2}}\right],     &   {r=0};\\
 &\\			
\left(\frac{r}{2^{l_1l_2}},~~\frac{r+1}{2^{l_1l_2}}\right],   &    {0<r \leq 2^{l_1l_2}-1}.
		\end{array} \right.
   \end{array}\right\}
 \end{eqnarray*}
 and
 \begin{eqnarray*}
 B_{r,\mathbf{l},\mathbf{\alpha},\mathbf{n}}=\{\omega:\varepsilon_{\mathbf{j}}(\omega)=
 \alpha_{\mathbf{j}},\mathbf{1}\leq \mathbf{j} \leq \mathbf{n}\}\bigcap B_{r,\mathbf{l}}.
 \end{eqnarray*}

\textbf{Lemma 4.1}. Under the conditions of Theorem 2.2, for $\mathbf{k,l} \in \mathbf{R_{\mathbf{n}}}$
such that $\mathbf{k}\neq \mathbf{l}$ and $l_1l_2\geq k_1k_2$, we have
\begin{eqnarray*}
&&\left|Cov \bigg(\mathbf{\mathbbm{1}}_{\{\bigcap_{\mathbf{i} \in \mathbf{R_{k}}}\{X_{\mathbf{i}}\leq u_{\mathbf{k,i}},X_{\mathbf{i}}(\varepsilon)\leq v_{\mathbf{k,i}}\}\}},\mathbf{\mathbbm{1}}_{\{\bigcap_{\mathbf{i} \in \mathbf{R_{l}}-\mathbf{R_{k}}}\{X_{\mathbf{i}}\leq u_{\mathbf{l,i}},X_{\mathbf{i}}(\varepsilon)\leq v_{\mathbf{l,i}}\}\}}\bigg)\right|\\
&&\ \ \ll\left\{
		\begin{array}{lcl}
			\alpha_{\mathbf{l,k},m_{l_1},m_{l_2}}^{*}+\frac{k_2m_{l_1}+k_1m_{l_2}+m_{l_1}m_{l_2}}{l_1l_2}
,     &     & {k_1< k_2,l_1<l_2};\\
& &\\			\alpha_{\mathbf{l,k},m_{l_1},m_{l_2}}^{*}+\frac{l_2m_{l_1}+m_{l_1}m_{l_2}}{l_1l_2},   &     & {k_1<l_1,l_2<k_2,k_1k_2<l_1l_2};\\
& &\\			\alpha_{\mathbf{l,k},m_{l_1},m_{l_2}}^{*}+\frac{l_1m_{l_2}+m_{l_1}m_{l_2}}{l_1l_2},   &     & {l_1<k_2,k_2<l_2,k_1k_2<l_1l_2}.
			\end{array} \right.
\end{eqnarray*}
\textbf{Proof. } Recall that $\mathbf{M_{kn}}=\{(j_1,j_2):(j_1,j_2)\in \mathbb{N}^2,0 \leq j_i \leq \#(\pi_i(\mathbf{M^{*}_{kn}}))+m_{n_i},i=1,2\}$ with $ \mathbf{M}_{\mathbf{kn}}^{*}
=\mathbf{R_k}\bigcap \mathbf{R_n}$. Write
\begin{eqnarray*}
&&\bigg|Cov \bigg(\mathbf{\mathbbm{1}}_{\{\bigcap_{\mathbf{i} \in \mathbf{R_{k}}}\{X_{\mathbf{i}}\leq u_{\mathbf{k,i}},X_{\mathbf{i}}(\varepsilon)\leq v_{\mathbf{k,i}}\}\}},\mathbf{\mathbbm{1}}_{\{\bigcap_{\mathbf{i} \in \mathbf{R_{l}}-\mathbf{R_{k}}}\{X_{\mathbf{i}}\leq u_{\mathbf{l,i}},X_{\mathbf{i}}(\varepsilon)\leq v_{\mathbf{l,i}}\}\}}\bigg)\bigg|\\
&&=\bigg|P\bigg(\mathcal{B}_\mathbf{k}(\mathbf{R_k})\bigcap\mathcal{B}_\mathbf{l}(\mathbf{R_l-R_k})\bigg)-P\bigg(\mathcal{B}_\mathbf{k}(\mathbf{R_k})\bigg)P\bigg(\mathcal{B}_\mathbf{l}(\mathbf{R_l-R_k})\bigg)\bigg|\\
&& \leq \bigg|P\bigg(\mathcal{B}_\mathbf{k}(\mathbf{R_k})\bigcap\mathcal{B}_\mathbf{l}(\mathbf{R_l-R_k})\bigg)-P\bigg(\mathcal{B}_\mathbf{k}(\mathbf{R_k})\bigcap\mathcal{B}_\mathbf{l}(\mathbf{R_l-M_{kl}})\bigg)\bigg|\\
&& \ \ +\bigg|P\bigg(\mathcal{B}_\mathbf{k}(\mathbf{R_k})\bigcap\mathcal{B}_\mathbf{l}(\mathbf{R_l-M_{kl}})\bigg)-P\bigg(\mathcal{B}_\mathbf{k}(\mathbf{R_k})\bigg)P\bigg(\mathcal{B}_\mathbf{l}(\mathbf{R_l-M_{kl}})\bigg)\bigg|\\
&& \ \ +\bigg|P\bigg(\mathcal{B}_\mathbf{k}(\mathbf{R_k})\bigg)P\bigg(\mathcal{B}_\mathbf{l}(\mathbf{R_l-M_{kl}})\bigg)-P\bigg(\mathcal{B}_\mathbf{k}(\mathbf{R_k})\bigg)P\bigg(\mathcal{B}_\mathbf{l}(\mathbf{R_l-R_k})\bigg)\bigg|\\
&& =I_1+I_2+I_3.
\end{eqnarray*}
Noting that $k_{1}k_{2}\leq l_{1}l_{2}$, we will estimate the above terms for three cases:
case 1, $k_1\leq l_1,k_2\leq l_2$; case 2, $k_1> l_1,k_2\leq l_2$ but $k_{1}k_{2}\leq l_{1}l_{2}$; case 3, $k_1\leq l_1,k_2> l_2$ but $k_{1}k_{2}\leq l_{1}l_{2}$.\\
For the first case, we have
\begin{eqnarray*}
\label{tan1}
I_1&=&\bigg|P\bigg(\mathcal{B}_\mathbf{k}(\mathbf{R_k})\bigcap\mathcal{B}_\mathbf{l}(\mathbf{R_l-R_k})\bigg)-P\bigg(\mathcal{B}_\mathbf{k}(\mathbf{R_k})\bigcap\mathcal{B}_\mathbf{l}(\mathbf{R_l-M_{kl}})\bigg)\bigg|\nonumber\\
&\leq& \bigg|P\bigg(\mathcal{B}_\mathbf{l}(\mathbf{R_l-R_k})\bigg)-P\bigg(\mathcal{B}_\mathbf{l}(\mathbf{R_l-M_{kl}})\bigg)\bigg|\nonumber\\
& \leq& P\left(\overline{\mathcal{B}}_\mathbf{l}(\mathbf{M_{kl}}-\mathbf{R_k})\right)\nonumber\\
& \leq& P\bigg(\bigcup_{\mathbf{i}\in \mathbf{M_{kl}}-\mathbf{R_k}}\left\{X_{\mathbf{i}}>u_{\mathbf{l,i}}\right\}\bigg)+
P\bigg(\bigcup_{\mathbf{i}\in \mathbf{M_{kl}}-\mathbf{R_k}}\left\{X_{\mathbf{i}}(\varepsilon)>v_{\mathbf{l,i}}\right\}\bigg)\nonumber\\
&=&I_{11}+I_{12},
\end{eqnarray*}
where
\begin{eqnarray*}
I_{11}&=&P\bigg(\bigcup_{\mathbf{i}\in \mathbf{M_{kl}}-\mathbf{R_k}}\left\{X_{\mathbf{i}}>u_{\mathbf{l,i}}\right\}\bigg)\leq
\sum_{\mathbf{i}\in \mathbf{M_{kl}}-\mathbf{R_k}}P\left(X_{\mathbf{i}}>u_{\mathbf{l,i}}\right)\\
&\leq& \left(k_2m_{l_1}+k_1m_{l_2}+m_{l_1}m_{l_2}\right)\max\left\{P\left(X_{\mathbf{i}}>u_{\mathbf{l,i}}\right),\mathbf{i}\leq \mathbf{l}\right\}\\
&=&\frac{k_2m_{l_1}+k_1m_{l_2}+m_{l_1}m_{l_2}}{l_1l_2} l_1l_2\max\left\{P\left(X_{\mathbf{i}}>v_{\mathbf{l,i}}\right),\mathbf{i}\leq \mathbf{l}\right\}\\
&\ll& \frac{k_2m_{l_1}+k_1m_{l_2}+m_{l_1}m_{l_2}}{l_1l_2}
\end{eqnarray*}
by using the condition that $\sup_{\mathbf{n}\geq\mathbf{1}}\{n_{1}n_{2}P(X_{\mathbf{i}}\geq v_{\mathbf{n},\mathbf{i}}), \mathbf{i}\leq \mathbf{n}\}$ is bounded, and similarly,
\begin{eqnarray*}
I_{12}&=&P\bigg(\bigcup_{\mathbf{i}\in \mathbf{M_{kl}}-\mathbf{R_k}}\left\{X_{\mathbf{i}}(\varepsilon)>v_{\mathbf{l,i}}\right\}\bigg)\\
&= &\sum_{t=0}^{\#(\mathbf{M_{kl}}-\mathbf{R_k})}P\left(\bigcup_{\mathbf{i}\in \mathbf{M_{kl}}-\mathbf{R_k}}X_{\mathbf{i}}(\varepsilon)>v_{\mathbf{l,i}}\bigg|\sum_{\mathbf{i}\in \mathbf{M_{kl}}-\mathbf{R_k}}\varepsilon_{\mathbf{i}}=t\right)P\left(\sum_{\mathbf{i}\in \mathbf{M_{kl}}-\mathbf{R_k}}\varepsilon_{\mathbf{i}}=t\right)\\
&\leq&\sum_{t=0}^{\#(\mathbf{M_{kl}}-\mathbf{R_k})}t\max\left\{P\left(X_{\mathbf{i}}>v_{\mathbf{l,i}}\right),\mathbf{i}\leq \mathbf{l}\right\}P\left(\sum_{\mathbf{i}\in \mathbf{M_{kl}}-\mathbf{R_k}}\varepsilon_{\mathbf{i}}=t\right)\\
& =&E\left(\sum_{\mathbf{i}\in \mathbf{M_{kl}}-\mathbf{R_k}}\varepsilon_{\mathbf{i}}\right)\max\left\{P\left(X_{\mathbf{i}}>v_{\mathbf{l,i}}\right),\mathbf{i}\leq \mathbf{l}\right\}\\
& =&\frac{E\left(\sum\limits_{\mathbf{i}\in \mathbf{M_{kl}}-\mathbf{R_k}}\varepsilon_{\mathbf{i}}\right)}{l_1l_2}l_1l_2\max\left\{P\left(X_{\mathbf{i}}>v_{\mathbf{l,i}}\right),\mathbf{i}\leq
\mathbf{l}\right\}\\
&\ll& \frac{k_2m_{l_1}+k_1m_{l_2}+m_{l_1}m_{l_2}}{l_1l_2}.
\end{eqnarray*}
Thus $$I_1\ll \frac{k_2m_{l_1}+k_1m_{l_2}+m_{l_1}m_{l_2}}{l_1l_2}.$$
By the similar arguments as for $I_{1}$, we have
$$I_3\ll \frac{k_2m_{l_1}+k_1m_{l_2}+m_{l_1}m_{l_2}}{l_1l_2}.$$
Let $\varepsilon(\mathbf{R_k})=\{\varepsilon_{\mathbf{i}}, \mathbf{i}\in\mathbf{R_k}\}$ and
$\varepsilon(\mathbf{R_l-M_{kl}})=\{\varepsilon_{\mathbf{i}}, \mathbf{i}\in\mathbf{R_l-M_{kl}}\}$.
By condition $\mathbf{D}^{*}\left(u_{\mathbf{k,j}}, v_{\mathbf{k,i}},u_{\mathbf{n,j}}, v_{\mathbf{n,i}}\right)$, we have
\begin{eqnarray}
\label{T1}
&&\bigg|E\left[P\bigg(\mathcal{B}_\mathbf{k}(\mathbf{R_k})\bigcap\mathcal{B}_\mathbf{l}(\mathbf{R_l-M_{kl}})|(\varepsilon(\mathbf{R_k}),\varepsilon(\mathbf{R_l-M_{kl}}))\bigg)\right]\nonumber\\
&&-E\left[P\bigg(\mathcal{B}_\mathbf{k}(\mathbf{R_k})|\varepsilon(\mathbf{R_k})\bigg)
P\bigg(\mathcal{B}_\mathbf{l}(\mathbf{R_l-M_{kl}})|\varepsilon(\mathbf{R_l-M_{kl}})\bigg)\right]\bigg|\leq \alpha^{*}_{\mathbf{l,k},m_{l_1},m_{l_2}}.
\end{eqnarray}
By the independence of $\{\varepsilon_{\mathbf{i}}, \mathbf{i\geq1}\}$, we have
\begin{eqnarray}
\label{T2}
&&E\left[P\bigg(\mathcal{B}_\mathbf{k}(\mathbf{R_k})|\varepsilon(\mathbf{R_k})\bigg)
P\bigg(\mathcal{B}_\mathbf{l}(\mathbf{R_l-M_{kl}})|\varepsilon(\mathbf{R_l-M_{kl}})\bigg)\right]\nonumber\\
&&= E\left[P\bigg(\mathcal{B}_\mathbf{k}(\mathbf{R_k})|\varepsilon(\mathbf{R_k})\bigg)\right]
E\left[P\bigg(\mathcal{B}_\mathbf{l}(\mathbf{R_l-M_{kl}})|\varepsilon(\mathbf{R_l-M_{kl}})\bigg)\right]\nonumber\\
&&=P\bigg(\mathcal{B}_\mathbf{k}(\mathbf{R_k})\bigg)
P\bigg(\mathcal{B}_\mathbf{l}(\mathbf{R_l-M_{kl}})\bigg).
\end{eqnarray}
Therefore (\ref{T1}) together with (\ref{T2}) implies
\begin{eqnarray*}
I_2=\bigg|P\bigg(\mathcal{B}_\mathbf{k}(\mathbf{R_k})\bigcap\mathcal{B}_\mathbf{l}(\mathbf{R_l-M_{kl}})\bigg)
-P\bigg(\mathcal{B}_\mathbf{k}(\mathbf{R_k})\bigg)P\bigg(\mathcal{B}_\mathbf{l}(\mathbf{R_l-M_{kl}})\bigg)\bigg|
\ll \alpha^{*}_{\mathbf{l,k},m_{l_1},m_{l_2}}.
\end{eqnarray*}
Thus, we have
\begin{eqnarray*}
&&\bigg|Cov \bigg(\mathbf{\mathbbm{1}}_{\{\bigcap_{\mathbf{i} \in \mathbf{R_{k}}}\{X_{\mathbf{i}}\leq u_{\mathbf{k,i}},X_{\mathbf{i}}(\varepsilon)\leq v_{\mathbf{k,i}}\}\}},\mathbf{\mathbbm{1}}_{\{\bigcap_{\mathbf{i} \in \mathbf{R_{l}}-\mathbf{R_{k}}}\{X_{\mathbf{i}}\leq u_{\mathbf{l,i}},X_{\mathbf{i}}(\varepsilon)\leq v_{\mathbf{l,i}}\}\}}\bigg)\bigg|\\
&&\ \ \ll \alpha^{*}_{\mathbf{l},\mathbf{k},m_{l_1},m_{l_2}}+\frac{k_1m_{l_2}+k_2m_{l_1}+m_{l_1}m_{l_2}}{l_1l_2}.
\end{eqnarray*}
Next, we deal with the second case: $k_1<l_1, l_2<k_2$, but $k_1k_2<l_1l_2$. As for $(\ref{tan1})$, we have
\begin{eqnarray*}
I_1&=&\bigg|P\bigg(\mathcal{B}_\mathbf{k}(\mathbf{R_k})\bigcap\mathcal{B}_\mathbf{l}(\mathbf{R_l-R_k})\bigg)-P\bigg(\mathcal{B}_\mathbf{k}(\mathbf{R_k})\bigcap\mathcal{B}_\mathbf{l}(\mathbf{R_l-M_{kl}})\bigg)\bigg|\\
& \leq& P\bigg(\bigcup_{\mathbf{i}\in \mathbf{M_{kl}}-\mathbf{R_k}}\left\{X_{\mathbf{i}}>u_{\mathbf{l,i}}\right\}\bigg)+
P\bigg(\bigcup_{\mathbf{i}\in \mathbf{M_{kl}}-\mathbf{R_k}}\left\{X_{\mathbf{i}}(\varepsilon)>v_{\mathbf{l,i}}\right\}\bigg)\\
&=&I_{13}+I_{14},
\end{eqnarray*}
where
\begin{eqnarray*}
I_{13}&\leq& (l_2m_{l_1}+m_{l_1}m_{l_2})\max\left\{P(X_{\mathbf{i}}>u_{\mathbf{l,i}}),\mathbf{i}\leq\mathbf{l}\right\}\\
&\leq&\frac{l_2m_{l_1}+m_{l_1}m_{l_2}}{l_1l_2}l_1l_2\max\left\{P(X_{\mathbf{i}}>v_{\mathbf{l,i}}),\mathbf{i}\leq\mathbf{l}\right\}\\
& \ll& \frac{l_2m_{l_1}+m_{l_1}m_{l_2}}{l_1l_2},
\end{eqnarray*}
by using the condition that $\sup_{\mathbf{n}\geq\mathbf{1}}\{n_{1}n_{2}P(X_{\mathbf{i}}\geq v_{\mathbf{n},\mathbf{i}}), \mathbf{i}\leq \mathbf{n}\}$ is bounded again, and similarly,
\begin{eqnarray*}
I_{14}&=& \sum_{t=0}^{\#(\mathbf{M_{kl}}-\mathbf{R_k})}P\left(\bigcup_{\mathbf{i}\in\mathbf{M_{kl}}-\mathbf{R_k}}X_{\mathbf{i}}(\varepsilon)>v_{\mathbf{l,i}}\bigg|\sum_{\mathbf{i}\in \mathbf{M_{kl}}-\mathbf{R_k}}\varepsilon_{\mathbf{i}}=t\right)P\left(\sum_{\mathbf{i}\in \mathbf{M_{kl}}-\mathbf{R_k}}\varepsilon_{\mathbf{i}}=t\right)\\
&\leq&\sum_{t=0}^{\#(\mathbf{M_{kl}}-\mathbf{R_k})}t\max\left\{P(X_{\mathbf{i}}>v_{\mathbf{l,i}}),\mathbf{i}\leq \mathbf{l}\right\}P\left(\sum_{\mathbf{i}\in \mathbf{M_{kl}}-\mathbf{R_k}}\varepsilon_{\mathbf{i}}=t\right)\\
& =&E\left(\sum_{\mathbf{i}\in \mathbf{M_{kl}}-\mathbf{R_k}}\varepsilon_{\mathbf{i}}\right)\max\left\{P(X_{\mathbf{i}}>v_{\mathbf{l,i}}),\mathbf{i}\leq \mathbf{l}\right\}\\
& =&\frac{E\left(\sum\limits_{\mathbf{i}\in \mathbf{M_{kl}}-\mathbf{R_k}}\varepsilon_{\mathbf{i}}\right)}{l_1l_2}l_1l_2\max\left\{P(X_{\mathbf{i}}>v_{\mathbf{l,i}}),\mathbf{i}\leq
\mathbf{l}\right\}\\
&\leq& \frac{l_2m_{l_1}+m_{l_1}m_{l_2}}{l_1l_2}l_1l_2\max\left\{P(X_{\mathbf{i}}>v_{\mathbf{l,i}}),\mathbf{i}\leq\mathbf{l}\right\}\\
&\ll& \frac{l_2m_{l_1}+m_{l_1}m_{l_2}}{l_1l_2}.
\end{eqnarray*}
Thus $$I_1\ll \frac{l_2m_{l_1}+m_{l_1}m_{l_2}}{l_1l_2}.$$
By the similar arguments as for $I_{1}$, we have
\begin{eqnarray*}
I_3 \ll \frac{l_2m_{l_1}+m_{l_1}m_{l_2}}{l_1l_2}.
\end{eqnarray*}
By condition $\mathbf{D}^{*}\left(u_{\mathbf{k,j}}, v_{\mathbf{k,i}},u_{\mathbf{n,j}}, v_{\mathbf{n,i}}\right)$ and the independence of $\{\varepsilon_{\mathbf{i}}, \mathbf{i\geq1}\}$ again, we get
\begin{eqnarray*}
I_2=\bigg|P\bigg(\mathcal{B}_\mathbf{k}(\mathbf{R_k})\bigcap\mathcal{B}_\mathbf{l}(\mathbf{R_l-M_{kl}})\bigg)-P\bigg(\mathcal{B}_\mathbf{k}(\mathbf{R_k})\bigg)P\bigg(\mathcal{B}_\mathbf{l}(\mathbf{R_l-M_{kl}})\bigg)\bigg|
\ll \alpha^{*}_{\mathbf{l,k},m_{l_1},m_{l_2}}.
\end{eqnarray*}
Thus, we have
\begin{eqnarray*}
&&\bigg|Cov \bigg(\mathbf{\mathbbm{1}}_{\{\bigcap_{\mathbf{i} \in \mathbf{R_{k}}}\{X_{\mathbf{i}}\leq u_{\mathbf{k,i}},X_{\mathbf{i}}(\varepsilon)\leq v_{\mathbf{k,i}}\}\}},\mathbf{\mathbbm{1}}_{\{\bigcap_{\mathbf{i} \in \mathbf{R_{l}}-\mathbf{R_{k}}}\{X_{\mathbf{i}}\leq u_{\mathbf{l,i}},X_{\mathbf{i}}(\varepsilon)\leq v_{\mathbf{l,i}}\}\}}\bigg)\bigg|\\
&&\ \ \ll \alpha^{*}_{\mathbf{l,k},m_{l_1},m_{l_2}}+\frac{l_2m_{l_1}+m_{l_1}m_{l_2}}{l_1l_2}.
\end{eqnarray*}
By the similar discussions as for the second case, we can get the desired bounds for the third case, so we omit the details. The proof of the lemma is complete.

\textbf{Lemma 4.2} Under the conditions of Theorem 2.3, for $\mathbf{k,l} \in \mathbf{R_{\mathbf{n}}}$
such that $\mathbf{k}\neq \mathbf{l}$ and $k_1k_2 \leq l_1l_2$, we have
\begin{equation}
E\left|\mathbf{\mathbbm{1}}_{\bigcap_{\mathbf{i}\in \mathbf{R_l}-\mathbf{R_k}}\left\{X_{\mathbf{i}}\leq u_{\mathbf{l,i}},X_{\mathbf{i}}(\varepsilon)\leq v_{\mathbf{l,i}}\right\}}-\mathbf{\mathbbm{1}}_{\bigcap_{\mathbf{i}\in \mathbf{R_l}}\left\{X_{\mathbf{i}}\leq u_{\mathbf{l,i}},X_{\mathbf{i}}(\varepsilon)\leq v_{\mathbf{l,i}}\right\}}\right|\leq \frac{l_1l_2-\#\left(\mathbf{R_l}-\mathbf{R_k}\right)}{l_1l_2}.
\end{equation}

\textbf{Proof.} We have
\begin{eqnarray*}
&&E\left|\mathbf{\mathbbm{1}}_{\bigcap_{\mathbf{i}\in \mathbf{R_l}-\mathbf{R_k}}\left\{X_{\mathbf{i}}\leq u_{\mathbf{l,i}},X_{\mathbf{i}}(\varepsilon)\leq v_{\mathbf{l,i}}\right\}}-\mathbf{\mathbbm{1}}_{\bigcap_{\mathbf{i}\in \mathbf{R_l}}\left\{X_{\mathbf{i}}\leq u_{\mathbf{l,i}},X_{\mathbf{i}}(\varepsilon)\leq v_{\mathbf{l,i}}\right\}}\right|\\
&&=P\left(\bigcap_{\mathbf{i}\in \mathbf{R_l}-\mathbf{R_k}}\left\{X_{\mathbf{i}}\leq u_{\mathbf{l,i}},X_{\mathbf{i}}(\varepsilon)\leq v_{\mathbf{l,i}}\right\}\right)-P\left(\bigcap_{\mathbf{i}\in \mathbf{R_l}}\left\{X_{\mathbf{i}}\leq u_{\mathbf{l,i}},X_{\mathbf{i}}(\varepsilon)\leq v_{\mathbf{l,i}}\right\}\right)\\
&&\leq \sum_{\mathbf{i}\in \mathbf{R_l}-\left(\mathbf{R_l}-\mathbf{R_k}\right)}P\left(\left\{X_{\mathbf{i}}>u_{\mathbf{l,i}}\right\}\bigcup
\left\{X_{\mathbf{i}}(\varepsilon)>v_{\mathbf{l,i}}\right\}\right)\\
&&\leq \left[l_1l_2-\#\left(\mathbf{R_l}-\mathbf{R_k}\right)\right]\left[\max
\left\{P\left(X_{\mathbf{i}}>u_{\mathbf{l,i}}\right),\mathbf{i}\leq \mathbf{l}\right\}+\max
\left\{P\left(X_{\mathbf{i}}(\varepsilon)>v_{\mathbf{l,i}}\right),\mathbf{i}\leq \mathbf{l}\right\}\right]\\
&&\leq 2\left[l_1l_2-\#\left(\mathbf{R_l}-\mathbf{R_k}\right)\right]\max
\left\{P\left(X_{\mathbf{i}}>v_{\mathbf{l,i}}\right),\mathbf{i}\leq \mathbf{l}\right\}\\
&&\ll \frac{l_1l_2-\#\left(\mathbf{R_l}-\mathbf{R_k}\right)}{l_1l_2},
\end{eqnarray*}
by using the condition that $\sup_{\mathbf{n}\geq\mathbf{1}}\{n_{1}n_{2}P(X_{\mathbf{i}}\geq v_{\mathbf{n},\mathbf{i}}), \mathbf{i}\leq \mathbf{n}\}$ is bounded.

The following lemma is from Tan and Wang (2013), which plays a crucial role in the proof of Theorem 2.2.

\textbf{Lemma 4.3}. Let $\eta_{\mathbf{i}},\mathbf{i}\in \mathbb{Z}_{+}^2$ be
uniformly bounded variables. Assume that
\begin{equation}
Var\left(\frac{1}{\log n_1\log n_2}\sum_{\mathbf{k}\in \mathbf{R_n}}\frac{1}{k_1k_2}\eta_{\mathbf{k}}\right)\ll \frac{1}{\left(\log\log n_1\log\log n_2\right)^{1+\epsilon}}
\nonumber
\end{equation}
Then
\begin{equation}
\frac{1}{\log n_1\log n_2}\sum_{\mathbf{k}\in \mathbf{R_n}}\frac{1}{k_1k_2}\left(\eta_{\mathbf{k}}-E\left(\eta_{\mathbf{k}}\right)\right)
\rightarrow 0    ~~~~~a.s.
\end{equation}
 \textbf{Proof}: Please see Lemma 3.2 of Tan and Wang (2013).

\textbf{Proof of Theorem 2.1:} Recall that $X_{\mathbf{i}}(\alpha)=(1-\alpha_{\mathbf{i}})\gamma(X_{\mathbf{i}})+\alpha_{\mathbf{i}}X_{\mathbf{i}}$.
Let $w(s_{1},s_{2})=\sharp(\mathbf{K_{s}})$.
It is easy to see that $w(s_{1},s_{2})\rightarrow\infty$ as $\mathbf{n}\rightarrow\infty$.
%We divide the interval $\{1,2,\ldots,n_1\}$ into $K_{s}, s=1,2,\ldots,k_{n_1}$ intervals of length $m_1=\lfloor n_1/k_{n_1}\rfloor$(the integer part of $n_1/k_{n_1})$, and divide the interval $\{1,2,\ldots,n_2\}$ into $K_{t}^{*}, t=1,2,\ldots,k_{n_2}$ intervals of length $m_2=\lfloor n_2/k_{n_2}\rfloor$, respectively. Thus we divide $\mathbf{R_n}$ into $\pi_1^{-1}(K_s)\bigcap \pi_2^{-1}(K_t^{*}), s, t=1,\ldots, k_{n_2}$ rectangles of area $m= m_1m_2$  as much as possible. We denote
%$$K_s=\{j:(s-1)m_1+1\leq j <sm_1\},~~1\leq s \leq k_{n_1},$$
%$$K_t^{*}=\{j:(t-1)m_2+1\leq j <tm_2\},~~1\leq t \leq k_{n_2}.$$
%and $S_{(s,t)}$ denotes the numbers of random variables that can be observed in the rectangular set $\pi_1^{-1}(K_1\bigcup K_2\bigcup \ldots K_s)\bigcap \pi_2^{-1}(K_1^*\bigcup K_2^*\bigcup \ldots K_t^*)$.\\
By using the full probability formula and triangle inequality, we get
\begin{eqnarray*}
  &&\left|P\left(\bigcap_{\mathbf{i}\in \mathbf{R_{n}}}\{X_{\mathbf{i}}(\varepsilon)\leq v_{\mathbf{n,i}}\},\bigcap_{\mathbf{i}\in \mathbf{R_{n}}}\{X_{\mathbf{i}}\leq u_{\mathbf{n,i}}\}\right)-E[e^{-\lambda\kappa}e^{-(1-\lambda)\tau}] \right|\\
  %&=&\sum_{r=0}^{2^{k_{n_1}k_{n_2}}-1}\sum_{\mathbf{\alpha} \in \{0,1\}^{n_1n_2}}\left|E\left(P\left(\bigcap_{\mathbf{i}\in \mathbf{R_{n}}}\{X_{\mathbf{i}}(\alpha)\leq v_{\mathbf{n,i}}\},\bigcap_{\mathbf{i}\in \mathbf{R_{n}}}\{X_{\mathbf{i}}\leq u_{\mathbf{n,i}}\}\right)-E[e^{-\lambda\kappa}e^{-(1-\lambda)\tau}]
 % \right)\right|\mathbbm{1}_{\{B_{r,k_{n_1}k_{n_2},\alpha,n}\}}\\
  &&\leq \sum_{r=0}^{2^{k_{n_1}k_{n_2}}-1}\sum_{\mathbf{\alpha} \in \{0,1\}^{n_1n_2}}E\left|P\left(\bigcap_{\mathbf{i}\in \mathbf{R_{n}}}\{X_{\mathbf{i}}(\alpha)\leq v_{\mathbf{n,i}}\},\bigcap_{\mathbf{i}\in \mathbf{R_{n}}}\{X_{\mathbf{i}}\leq u_{\mathbf{n,i}}\}\right)\right.\\
  &&\ \ \ \left.-\prod \limits_{s_1=1}^{k_{n_1}}\prod \limits_{s_2=1}^{k_{n_2}}
  P\left(\bigcap_{\mathbf{i}\in \mathbf{K_{s}}}\{X_{\mathbf{i}}(\alpha)\leq v_{\mathbf{n,i}}\},\bigcap_{\mathbf{i}\in \mathbf{K_s}}\{X_{\mathbf{i}}\leq u_{\mathbf{n,i}}\}\right)\right|\mathbbm{1}_{\{B_{r,\mathbf{k_{n}},\alpha,\mathbf{n}}\}}  \\
  &&+\sum_{r=0}^{2^{k_{n_1}k_{n_2}}-1}\sum_{\mathbf{\alpha} \in \{0,1\}^{n_1n_2}}E\left|\prod \limits_{s_1=1}^{k_{n_1}}\prod \limits_{s_2=0}^{k_{n_2}}
  P\left(\bigcap_{\mathbf{i}\in \mathbf{K_{s}}}\{X_{\mathbf{i}}(\alpha)\leq v_{\mathbf{n,i}}\},\bigcap_{\mathbf{i}\in \mathbf{K_{s}}}\{X_{\mathbf{i}}\leq u_{\mathbf{n,i}}\}\right)\right.\\
  &&\ \ \ -\left.\prod \limits_{s_1=1}^{k_{n_1}}\prod \limits_{s_2=1}^{k_{n_2}} \left[1-\frac{\frac{r}{2^{k_{n_1}k_{n_2}}}\sum\limits_{\mathbf{i}
  \in \mathbf{R}_{\mathbf{n}}}
  P(X_{\mathbf{i}}>v_{\mathbf{n,i}})+\left(1-\frac{r}
  {2^{k_{n_1}k_{n_2}}}\right)\sum\limits_{\mathbf{i}
  \in \mathbf{R}_{\mathbf{n}}}P(X_{\mathbf{i}}>u_{\mathbf{n,i}})}{k_{n_1}k_{n_2}}
  \right]\right|\mathbbm{1}_{\{B_{r,\mathbf{k_{n}},\alpha,\mathbf{n}}\}}\\
  &&+\sum_{r=0}^{2^{k_{n_1}k_{n_2}}-1}\sum_{\mathbf{\alpha} \in \{0,1\}^{n_1n_2}}E\left|\prod \limits_{s_1=1}^{k_{n_1}}\prod \limits_{s_2=1}^{k_{n_2}}\left[1-\frac{\frac{r}{2^{k_{n_1}k_{n_2}}}
  \sum\limits_{\mathbf{i}
  \in \mathbf{R}_{\mathbf{n}}}P(X_{\mathbf{i}}>v_{\mathbf{n,i}})+\left(1-\frac{r}
  {2^{k_{n_1}k_{n_2}}}\right)\sum\limits_{\mathbf{i}
  \in \mathbf{R}_{\mathbf{n}}}P(X_{\mathbf{i}}>u_{\mathbf{n,i}})}{k_{n_1}k_{n_2}}
  \right]\right.\\
  &&\ \ \ -\left.\prod \limits_{s_1=1}^{k_{n_1}}\prod \limits_{s_2=1}^{k_{n_2}} \left[1-\frac{\lambda \sum\limits_{\mathbf{i}
  \in \mathbf{R}_{\mathbf{n}}}
 P(X_{\mathbf{i}}>v_{\mathbf{n,i}})+\left(1-\lambda\right)
 \sum\limits_{\mathbf{i}
  \in \mathbf{R}_{\mathbf{n}}}
 P(X_{\mathbf{i}}>u_{\mathbf{n,i}})}{k_{n_1}k_{n_2}}
  \right]\right|\mathbbm{1}_{\{B_{r,\mathbf{k_{n}},\alpha,\mathbf{n}}\}}\\
  &&+\sum_{r=0}^{2^{k_{n_1}k_{n_2}}-1}\sum_{\mathbf{\alpha} \in \{0,1\}^{n_1n_2}}E\left|\prod \limits_{s_1=1}^{k_{n_1}}\prod \limits_{s_2=1}^{k_{n_2}}\left[1-\frac{\lambda \sum\limits_{\mathbf{i}
  \in \mathbf{R}_{\mathbf{n}}}
  P(X_{\mathbf{i}}>v_{\mathbf{n,i}})+\left(1-\lambda\right)\sum\limits
  _{\mathbf{i}
  \in \mathbf{R}_{\mathbf{n}}}
 P(X_{\mathbf{i}}>u_{\mathbf{n,i}})}{k_{n_1}k_{n_2}}
  \right]\right.\\
  &&\ \ \ \left.-E[e^{-\lambda\kappa}e^{-(1-\lambda)\tau}]\right|
  \mathbbm{1}_{\{B_{r,\mathbf{k_{n}},\alpha,\mathbf{n}}\}}\\
  &&=J_1+J_2+J_3+J_4.
  \end{eqnarray*}
  To bound the first term  $J_1$, we will divide each rectangle subsets $\mathbf{K_{s}}=\mathbf{K}_{(s_{1},s_{2})}$, $s_{1}=1,2,\ldots, k_{n_{1}}$, $s_{2}=1,2,\ldots, k_{n_{2}}$ into two parts. Without loss of generality, suppose the coordinates of the four apexes of the rectangle $\mathbf{K_{s}}$ are $(s_{11},s_{21}),(s_{11},s_{22}),(s_{12},s_{21})$ and $(s_{12},s_{22})$ with $s_{11}<s_{12}$ and $s_{21}<s_{22}$. Let $\mathbf{K_{s}^*}=[s_{11},s_{12}-m_{n_1}]\times[s_{21},s_{22}-m_{n_2}]$ and $\mathbf{K_{s}^{**}}=\mathbf{K_{s}}-\mathbf{K_{s}^*}$.
 It is easy to see that
 \begin{eqnarray*}
 \sharp\left(\bigcup_{s_{1}=1,2,\ldots, k_{n_{1}},s_{1}=1,2,\ldots, k_{n_{1}}}\mathbf{K_{s}^{**}}\right)<m_{n_1}k_{n_{1}}n_2+m_{n_2}k_{n_{2}}n_1.
   \end{eqnarray*}
Obviously,
 \begin{eqnarray*}
&&\left|P\left(\bigcap_{\mathbf{i}\in \mathbf{R_{n}}}\{X_{\mathbf{i}}(\alpha)\leq v_{\mathbf{n,i}}\},\bigcap_{\mathbf{i}\in \mathbf{R_{n}}}\{X_{\mathbf{i}}\leq u_{\mathbf{n,i}}\}\right)-\prod \limits_{s_1=1}^{k_{n_1}}\prod \limits_{s_2=1}^{k_{n_2}}
  P\left(\bigcap_{\mathbf{i}\in \mathbf{K_{s}}}\{X_{\mathbf{i}}(\alpha)\leq v_{\mathbf{n,i}}\},\bigcap_{\mathbf{i}\in \mathbf{K_s}}\{X_{\mathbf{i}}\leq u_{\mathbf{n,i}}\}\right)\right|\\
  &&\leq \left|P\left(\bigcap_{\mathbf{i}\in \mathbf{R_{n}}}\{X_{\mathbf{i}}(\alpha)\leq v_{\mathbf{n,i}}\},\bigcap_{\mathbf{i}\in \mathbf{R_{n}}}\{X_{\mathbf{i}}\leq u_{\mathbf{n,i}}\}\right)-P\left(\bigcap_{\mathbf{i}\in \cup \mathbf{K_{s}^{*}}}\{X_{\mathbf{i}}(\alpha)\leq v_{\mathbf{n,i}}\},\bigcap_{\mathbf{i}\in \cup \mathbf{K_{s}^{*}}}\{X_{\mathbf{i}}\leq u_{\mathbf{n,i}}\}\right)\right|\\
  &&+\left|P\left(\bigcap_{\mathbf{i}\in \cup \mathbf{K_{s}^{*}}}\{X_{\mathbf{i}}(\alpha)\leq v_{\mathbf{n,i}}\},\bigcap_{\mathbf{i}\in \cup \mathbf{K_{s}^{*}}}\{X_{\mathbf{i}}\leq u_{\mathbf{n,i}}\}\right)-\prod \limits_{s_1=1}^{k_{n_1}}\prod \limits_{s_2=1}^{k_{n_2}}
  P\left(\bigcap_{\mathbf{i}\in \mathbf{K_{s}^{*}}}\{X_{\mathbf{i}}(\alpha)\leq v_{\mathbf{n,i}}\},\bigcap_{\mathbf{i}\in \mathbf{K_s}^{*}}\{X_{\mathbf{i}}\leq u_{\mathbf{n,i}}\}\right)\right|\\
  &&+\left|\prod \limits_{s_1=1}^{k_{n_1}}\prod \limits_{s_2=1}^{k_{n_2}}
  P\left(\bigcap_{\mathbf{i}\in \mathbf{K_{s}^{*}}}\{X_{\mathbf{i}}(\alpha)\leq v_{\mathbf{n,i}}\},\bigcap_{\mathbf{i}\in \mathbf{K_s}^{*}}\{X_{\mathbf{i}}\leq u_{\mathbf{n,i}}\}\right)-\prod \limits_{s_1=1}^{k_{n_1}}\prod \limits_{s_2=1}^{k_{n_2}}
  P\left(\bigcap_{\mathbf{i}\in \mathbf{K_{s}}}\{X_{\mathbf{i}}(\alpha)\leq v_{\mathbf{n,i}}\},\bigcap_{\mathbf{i}\in \mathbf{K_s}}\{X_{\mathbf{i}}\leq u_{\mathbf{n,i}}\}\right)\right|\\
  &&=J_{11}+J_{12}+J_{13}.
\end{eqnarray*}
It is directly to check that
 \begin{eqnarray*}
 J_{11}&\leq& P\left(\bigcup_{\mathbf{i}\in \cup \mathbf{K_{s}^{**}}}\{X_{\mathbf{i}}(\alpha)> v_{\mathbf{n,i}}\}\right)+P\left(\bigcup_{\mathbf{i}\in \cup \mathbf{K_{s}^{**}}}\{X_{\mathbf{i}}> u_{\mathbf{n,i}}\}\right)\\
 &\leq&2 \sharp\left(\cup\mathbf{K_{s}^{**}}\right)\max\limits_{\mathbf{i}\in \mathbf{R}_{\mathbf{n}}}P(X_{\mathbf{i}}>v_{\mathbf{n,i}})
 < 2[m_{n_1}k_{n_{1}}n_2+m_{n_2}k_{n_{2}}n_1] \max\limits_{\mathbf{i}\in \mathbf{R}_{\mathbf{n}}}P(X_{\mathbf{i}}>v_{\mathbf{n,i}}).
 \end{eqnarray*}
  Noting that
  \begin{eqnarray}
  \label{inequality}
  \left|\prod\limits_{s=1}^{k}a_s-\prod\limits_{s=1}^{k}b_s\right| \leq \sum\limits_{s=1}^{k}|a_s-b_s|,
  \end{eqnarray}
  for all $a_s,b_s \in [0,1]$, we have similarly
  \begin{eqnarray*}
 J_{13}&\leq&  \sum_{s_1=1}^{k_{n_1}}\sum_{s_2=1}^{k_{n_2}} \left[P\left(\bigcup_{\mathbf{i}\in \mathbf{K_{s}^{**}}}\{X_{\mathbf{i}}(\alpha)> v_{\mathbf{n,i}}\}\right)+P\left(\bigcup_{\mathbf{i}\in \mathbf{K_{s}^{**}}}\{X_{\mathbf{i}}> u_{\mathbf{n,i}}\}\right)\right]\\
 &\leq&2 \sharp\left(\cup\mathbf{K_{s}^{**}}\right)\max\limits_{\mathbf{i}\in \mathbf{R}_{\mathbf{n}}}P(X_{\mathbf{i}}>v_{\mathbf{n,i}})
 < 2[m_{n_1}k_{n_{1}}n_2+m_{n_2}k_{n_{2}}n_1] \max\limits_{\mathbf{i}\in \mathbf{R}_{\mathbf{n}}}P(X_{\mathbf{i}}>v_{\mathbf{n,i}})
 \end{eqnarray*}
 By induction and the condition $\mathbf{D}(u_{\mathbf{n,i}},v_{\mathbf{n,i}})$, we get
 \begin{eqnarray*}
 J_{12}&\leq&(k_{n_1}k_{n_2}-1)\alpha_{\mathbf{n},m_{n_1},m_{n_2}}.
 \end{eqnarray*}
 Thus, we have
  \begin{eqnarray*}
  J_1 &\leq& (k_{n_1}k_{n_2}-1)\alpha_{\mathbf{n},m_{n_1},m_{n_2}}+
 \frac{4[k_{n_1}m_{n_1}n_2+k_{n_2}m_{n_2}n_1]}{n_1n_2}
  n_1n_2\max\limits_{\mathbf{i}\in \mathbf{R}_{\mathbf{n}}}P(X_{\mathbf{i}}>v_{\mathbf{n,i}}).
  \end{eqnarray*}
  Further, noting that $\sup_{\mathbf{n}\geq\mathbf{1}}\{n_{1}n_{2}P(X_{\mathbf{i}}\geq v_{\mathbf{n},\mathbf{i}}), \mathbf{i}\leq \mathbf{n}\}$ is bounded and using condition $\mathbf{D}(u_{\mathbf{n,i}},v_{\mathbf{n,i}})$, we have
  $$J_1=o(1)~~\mbox{as}~~~\mathbf{n}\rightarrow \mathbf{\infty}.$$
  For the second term, for any $0 \leq r \leq 2^{k_{n_1}k_{n_2}}-1$, it is not hard to check from (\ref{eqTan1}) and Bonferroni type inequality that
  \begin{eqnarray*}
  &&\left[1-\frac{\frac{r}{2^{k_{n_1}k_{n_2}}}\sum\limits_{\mathbf{i}\in
   \mathbf{R}_{\mathbf{n}}}P(X_{\mathbf{i}}>v_{\mathbf{n,i}})
+\left(1-\frac{r}{2^{k_{n_1}k_{n_2}}}\right)\sum\limits_{\mathbf{i}\in
   \mathbf{R}_{\mathbf{n}}}P(X_{\mathbf{i}}>u_{\mathbf{n,i}})}{k_{n_1}k_{n_2}}
   \right]\\
   &&\ \ +\sum_{\mathbf{i}\in \mathbf{K_{s}}}\left[\frac{r}{2^{k_{n_1}k_{n_2}}}-\alpha_{\mathbf{i}}\right]\left(P(X_{\mathbf{i}}>v_{\mathbf{n,i}})-
   P(X_{\mathbf{i}}>u_{\mathbf{n,i}})\right)+o(1)\\
   &&\leq P\left(\bigcap_{\mathbf{i}\in \mathbf{K_{s}}}\{X_{\mathbf{i}}(\alpha)\leq v_{\mathbf{n,i}}\},\bigcap_{\mathbf{i}\in \mathbf{K_{s}}}\{X_{\mathbf{i}}\leq u_{\mathbf{n,i}}\}\right)\\
    &&\leq\left[1-\frac{\frac{r}{2^{k_{n_1}k_{n_2}}}\sum\limits_{\mathbf{i}\in
   \mathbf{R}_{\mathbf{n}}}P(X_{\mathbf{i}}>v_{\mathbf{n,i}})
+\left(1-\frac{r}{2^{k_{n_1}k_{n_2}}}\right)\sum\limits_{\mathbf{i}\in
   \mathbf{R}_{\mathbf{n}}}P(X_{\mathbf{i}}>u_{\mathbf{n,i}})}{k_{n_1}k_{n_2}}
   \right]\\
   &&\ \ +\sum_{\mathbf{i}\in \mathbf{K_{s}}}\left[\frac{r}{2^{k_{n_1}k_{n_2}}}-\alpha_{\mathbf{i}}\right]\left(P(X_{\mathbf{i}}>v_{\mathbf{n,i}})-
   P(X_{\mathbf{i}}>u_{\mathbf{n,i}})\right)\\
   &&\ \ +\sum\limits_{\mathbf{i,j} \in \mathbf{K_{s}},\mathbf{i}\neq \mathbf{j}}P\left(X_{\mathbf{i}}\geq v_{\mathbf{n,i}},X_{\mathbf{j}}\geq v_{\mathbf{n,j}}\right)+o(1),
    \end{eqnarray*}
  as $\mathbf{n}\rightarrow\infty$, which together with (\ref{inequality}) implies
  \begin{eqnarray*}
  J_2 &\leq & \sum\limits_{r=0}^{2^{k_{n_1}k_{n_2}}-1}\sum\limits_{\mathbf{\alpha} \in \{0,1\}^{n_1n_2}}\sum\limits_{s_1=1}^{k_{n_1}}\sum\limits_{s_2=1}^{k_{n_2}}
  E\left|P\left(\bigcap_{\mathbf{i}\in \mathbf{K_{s}}}\{X_{\mathbf{i}}(\alpha)\leq v_{\mathbf{n,i}}\},\bigcap_{\mathbf{i}\in \mathbf{K_{s}}}\{X_{\mathbf{i}}\leq u_{\mathbf{n,i}}\}\right)\right.\\
  &&\left. -\left[1-\frac{\frac{r}{2^{k_{n_1}k_{n_2}}}\sum\limits_{\mathbf{i}\in
   \mathbf{R}_{\mathbf{n}}}P(X_{\mathbf{i}}>v_{\mathbf{n,i}})
+\left(1-\frac{r}{2^{k_{n_1}k_{n_2}}}\right)\sum\limits_{\mathbf{i}\in
   \mathbf{R}_{\mathbf{n}}}P(X_{\mathbf{i}}>u_{\mathbf{n,i}})}{k_{n_1}k_{n_2}}
   \right]\right|\mathbbm{1}_{\{B_{r,\mathbf{k_{n}},\alpha,\mathbf{n}}\}}\\
  &\leq&\sum\limits_{r=0}^{2^{k_{n_1}k_{n_2}}-1}\sum\limits_{\mathbf{\alpha} \in \{0,1\}^{n_1n_2}}\sum\limits_{s_1=1}^{k_{n_1}}\sum\limits_{s_2=1}^{k_{n_2}}
  E\left|\sum_{\mathbf{i}\in \mathbf{K_{s}}}\left[\frac{r}{2^{k_{n_1}k_{n_2}}}-\alpha_{\mathbf{i}}\right]\left(P(X_{\mathbf{i}}>v_{\mathbf{n,i}})-
   P(X_{\mathbf{i}}>u_{\mathbf{n,i}})\right)\right|\mathbbm{1}_{\{B_{r,\mathbf{k_{n}},\alpha,\mathbf{n}}\}}\\
   &&+\sum\limits_{r=0}^{2^{k_{n_1}k_{n_2}}-1}\sum\limits_{\mathbf{\alpha} \in \{0,1\}^{n_1n_2}}\sum\limits_{s_1=1}^{k_{n_1}}\sum\limits_{s_2=1}^{k_{n_2}}
   E\left(\sum\limits_{\mathbf{i,j} \in \mathbf{K_{s}},\mathbf{i}\neq \mathbf{j}}P\left(X_{\mathbf{i}}\geq v_{\mathbf{n,i}},X_{\mathbf{j}}\geq v_{\mathbf{n,j}}\right)\right)
   \mathbbm{1}_{\{B_{r,\mathbf{k_{n}},\alpha,\mathbf{n}}\}}+o(1)\\
   &=&\sum\limits_{r=0}^{2^{k_{n_1}k_{n_2}}-1}\sum\limits_{s_1=1}^{k_{n_1}}
   \sum\limits_{s_2=1}^{k_{n_2}}
   E\left|\sum_{\mathbf{i}\in \mathbf{K_{s}}}\left[\frac{r}{2^{k_{n_1}k_{n_2}}}-\varepsilon_{\mathbf{i}}\right]\left(P(X_{\mathbf{i}}>v_{\mathbf{n,i}})-
   P(X_{\mathbf{i}}>u_{\mathbf{n,i}})\right)\right|\mathbbm{1}_{\{B_{r,\mathbf{k_{n}}}\}}\\
   &&+\sum\limits_{s_1=1}^{k_{n_1}}\sum\limits_{s_2=1}^{k_{n_2}}\sum\limits_
   {\mathbf{i,j}\in \mathbf{K_{s}},\mathbf{i}\neq \mathbf{j}}P\left(X_{\mathbf{i}}\geq v_{\mathbf{n,i}},X_{\mathbf{j}}\geq v_{\mathbf{n,j}}\right)+o(1)\\
   &=&J_{21}+J_{22}+o(1).
  \end{eqnarray*}
  For the first term $J_{21}$,  suppose again that the coordinates of the four apexes of the rectangle $\mathbf{K_{s}}$ are $(s_{11},s_{21}),(s_{11},s_{22}),(s_{12},s_{21})$ and $(s_{12},s_{22})$ with $s_{11}<s_{12}$ and $s_{21}<s_{22}$. It follows from the facts
  that
 $P(X_{\mathbf{i}}\leq u_{\mathbf{n,i}})=O(P(X_{\mathbf{j}}\leq u_{\mathbf{n,j}}))$ and $P(X_{\mathbf{i}}\leq v_{\mathbf{n,i}})=O(P(X_{\mathbf{j}}\leq v_{\mathbf{n,j}}))$ for all $\mathbf{1\leq i\neq j\leq n}$ that
  \begin{eqnarray*}
  &&\sum\limits_{r=0}^{2^{k_{n_1}k_{n_2}}-1}
  E\left|\sum_{\mathbf{i}\in \mathbf{K_{s}}}\left[\frac{r}{2^{k_{n_1}k_{n_2}}}-\varepsilon_{\mathbf{i}}\right]\left(P(X_{\mathbf{i}}>v_{\mathbf{n,i}})-
   P(X_{\mathbf{i}}>u_{\mathbf{n,i}})\right)\right|
   \mathbbm{1}_{\{B_{r,\mathbf{k_{n}}}\}}\\
  &&\ll  E \left|\sum_{\mathbf{i} \in \mathbf{K_{s}}} \left[\varepsilon_{\mathbf{i}}-\lambda\right]\right|(P(X_{\mathbf{1}}\leq u_{\mathbf{n,1}})-P(X_{\mathbf{1}}\leq v_{\mathbf{n,1}}))\\
    &&\ \ \ +\sum_{r=0}^{2^{k_{n_1}k_{n_2}}-1}(P(X_{\mathbf{1}}\leq u_{\mathbf{n,1}})-P(X_{\mathbf{1}}\leq v_{\mathbf{n,1}}))E \left|\sum_{\mathbf{i} \in \mathbf{K_{s}}} \left[\lambda-\frac{r}{2^{k_{n_1}k_{n_2}}}\right]\right|\mathbbm{1}_{\{B_{r,\mathbf{k_{n}}}\}} \nonumber\\
&&\leq \frac{1}{k_{n_{1}}k_{n_{2}}}E\left|\sum_{\mathbf{i} \in \mathbf{K_{s}}} \frac{\varepsilon_{\mathbf{i}}}{w(s_{1},s_{2})}-\lambda\right|n_{1}n_{2}(P(X_{\mathbf{1}}\leq u_{\mathbf{n,1}})-P(X_{\mathbf{1}}\leq v_{\mathbf{n,1}}))+\frac{1}{2^{k_{n_1}k_{n_2}}}\frac{n_{1}n_{2}(P(X_{\mathbf{1}}\leq u_{\mathbf{n,1}})-P(X_{\mathbf{1}}\leq v_{\mathbf{n,1}}))}{k_{n_1}k_{n_2}},\nonumber
 \end{eqnarray*}
where
 \begin{eqnarray*}
E\left|\sum_{\mathbf{i} \in \mathbf{K_{s}}} \frac{\varepsilon_{\mathbf{i}}}{w(s_{1},s_{2})}-\lambda\right|
&=& E\left| \frac{S_{(s_{12},s_{22})}-S_{(s_{11},s_{22})}-S_{(s_{12},s_{21})}+S_{(s_{11},s_{21})}}
{w(s_{1},s_{2})\emph{}} - \lambda \right|\\
&\leq& E\left|s_{12}s_{22}\left(\frac{S(s_{12},s_{22})}{s_{12}s_{22}w(s_{1},s_{2})}-\lambda\right)\right|
+E\left|s_{11}s_{22}\left(\frac{S_{(s_{11},s_{22})}}{s_{11}s_{22}w(s_{1},s_{2})}-\lambda\right)\right|\\
& +&E\left|s_{12}s_{21}\left(\frac{S_{(s_{12},s_{21})}}{s_{12}s_{21}w(s_{1},s_{2})}-\lambda\right)\right|
+E\left|s_{11}s_{21}\left(\frac{S_{(s_{11},s_{21})}}{s_{11}s_{21}w(s_{1},s_{2})}-\lambda\right)\right|.
 \end{eqnarray*}
Taking into account (\ref{eqT11.5}), we have
\begin{eqnarray*}
&&\lim_{w(s_{1},s_{2})\rightarrow\infty}\left[E\left|s_{12}s_{22}\left(\frac{S(s_{12},s_{22})}{s_{12}s_{22}w(s_{1},s_{2})}-\lambda\right)\right|
+E\left|s_{11}s_{22}\left(\frac{S_{(s_{11},s_{22})}}{s_{11}s_{22}w(s_{1},s_{2})}-\lambda\right)\right|\right.\\
&&\ \ \ \ \ \ \ \ \ \ \ \  \left.+E\left|s_{12}s_{21}\left(\frac{S_{(s_{12},s_{21})}}{s_{12}s_{21}w(s_{1},s_{2})}-\lambda\right)\right|
+E\left|s_{11}s_{21}\left(\frac{S_{(s_{11},s_{21})}}{s_{11}s_{21}w(s_{1},s_{2})}-\lambda\right)\right|\right]=0
\end{eqnarray*}
Note that $w(s_{1},s_{2})\rightarrow\infty$ as $\mathbf{n}\rightarrow\infty$.
Letting $\mathbf{n}\rightarrow\infty$, we have
\begin{eqnarray*}
J_{21}&\leq& \sum\limits_{s_1=1}^{k_{n_1}}\sum\limits_{s_2=1}^{k_{n_2}}\frac{1}{2^{k_{n_1}k_{n_2}}}
\frac{\kappa-\tau}{k_{n_1}k_{n_2}}=O(\frac{1}{2^{k_{n_1}k_{n_2}}}).
\end{eqnarray*}
  For the second term $J_{22}$,  by the condition $\mathbf{D}'(u_{\mathbf{n,i}})$, we  obtain
  \begin{eqnarray*}
  J_{22}&=&\frac{1}{k_{n_1}
  k_{n_2}}\sum\limits_{s_1=1}^{k_{n_1}}\sum\limits_{s_2=1}^{k_{n_2}}\left(k_{n_1}
  k_{n_2}\sum\limits_
   {\mathbf{i,j}\in \mathbf{K_{s}},\mathbf{i}\neq \mathbf{j}}P\left(X_{\mathbf{i}}\geq v_{\mathbf{n,i}},X_{\mathbf{j}}\geq v_{\mathbf{n,j}}\right)\right)
   =o(1),
   \end{eqnarray*}
   as $\mathbf{n}\rightarrow\infty$.
  We thus have  $$J_2\leq O(\frac{1}{2^{k_{n_1}k_{n_2}}})~~\mbox{as}\ \ \mathbf{n}\rightarrow\infty.$$
  For the term $ J_3$, applying the facts that $\sum_{\mathbf{i}\in \mathbf{R_{n}}}P\left(X_\mathbf{i}>u_{\mathbf{n,i}}\right)\rightarrow\tau>0$ and $\sum_{\mathbf{i}\in \mathbf{R_{n}}}P\left(X_\mathbf{i}>v_{\mathbf{n,i}}\right)\rightarrow\kappa>0$ again, we have
  \begin{eqnarray*}
  J_3 &\leq &\sum_{r=0}^{2^{k_{n_1}k_{n_2}}-1}\sum_{\mathbf{\alpha} \in \{0,1\}^{n_1n_2}}\sum_{s_1=1}^{k_{n_1}}\sum_{s_2=1}^{k_{n_2}}E\left|
  \left[1-\frac{\frac{r}{2^{k_{n_1}k_{n_2}}}
  \sum\limits_{\mathbf{i}
  \in \mathbf{R}_{\mathbf{n}}}P(X_{\mathbf{i}}>v_{\mathbf{n,i}})+\left(1-\frac{r}
  {2^{k_{n_1}k_{n_2}}}\right)\sum\limits_{\mathbf{i}
  \in \mathbf{R}_{\mathbf{n}}}P(X_{\mathbf{i}}>u_{\mathbf{n,i}})}{k_{n_1}k_{n_2}}
  \right]\right.\\
  &-&\left. \left[1-\frac{\lambda \sum\limits_{\mathbf{i}
  \in \mathbf{R}_{\mathbf{n}}}
 P(X_{\mathbf{i}}>v_{\mathbf{n,i}})+\left(1-\lambda\right)
 \sum\limits_{\mathbf{i}
  \in \mathbf{R}_{\mathbf{n}}}
 P(X_{\mathbf{i}}>u_{\mathbf{n,i}})}{k_{n_1}k_{n_2}}
  \right]\right|\mathbbm{1}_{\{B_{r,\mathbf{k_{n}},\alpha,\mathbf{n}}\}}\\
  &\leq&\sum_{r=0}^{2^{k_{n_1}k_{n_2}}-1}\sum_{\mathbf{\alpha} \in \{0,1\}^{n_1n_2}}\sum_{s_1=1}^{k_{n_1}}\sum_{s_2=1}^{k_{n_2}}E\left|
  \lambda-\frac{r}{2^{k_{n_1}k_{n_2}}}\right|\mathbbm{1}_{\{B_{r,\mathbf{k_{n}},\alpha,\mathbf{n}}\}}
  \frac{\sum\limits_{\mathbf{i}
  \in \mathbf{R}_{\mathbf{n}}}\left(P(X_{\mathbf{i}}>v_{\mathbf{n,i}})+
  P(X_{\mathbf{i}}>u_{\mathbf{n,i}})\right)}{k_{n_1}k_{n_2}}\\
  &\leq&\frac{\sum\limits_{\mathbf{i}
  \in \mathbf{R}_{\mathbf{n}}}\left(P(X_{\mathbf{i}}>v_{\mathbf{n,i}})+
  P(X_{\mathbf{i}}>u_{\mathbf{n,i}})\right)}{2^{k_{n_1}k_{n_2}}}\\
  &\leq& O(\frac{1}{2^{k_{n_1}k_{n_2}}}),
  \end{eqnarray*}
  as $\mathbf{n}\rightarrow \mathbf{\infty}$.
  For the last term $J_4$, letting $\mathbf{n}\rightarrow\infty$ and then $(k_{n_1},k_{n_2})\rightarrow  \mathbf{\infty}$, we have
  \begin{eqnarray*}
  &&\prod \limits_{s_1=1}^{k_{n_1}}\prod \limits_{s_2=1}^{k_{n_2}}E\left[1-\frac{\lambda \sum\limits_{\mathbf{i}
  \in \mathbf{R}_{\mathbf{n}}}
  P(X_{\mathbf{i}}>v_{\mathbf{n,i}})+\left(1-\lambda\right)\sum\limits
  _{\mathbf{i}
  \in \mathbf{R}_{\mathbf{n}}}
 P(X_{\mathbf{i}}>u_{\mathbf{n,i}})}{k_{n_1}k_{n_2}}\right]\\
 &&=E\left[1-\frac{\lambda \sum\limits_{\mathbf{i}
  \in \mathbf{R}_{\mathbf{n}}}
  P(X_{\mathbf{i}}>v_{\mathbf{n,i}})+\left(1-\lambda\right)\sum\limits
  _{\mathbf{i}
  \in \mathbf{R}_{\mathbf{n}}}
 P(X_{\mathbf{i}}>u_{\mathbf{n,i}})}{k_{n_1}k_{n_2}}
  \right]^{k_{n_1}k_{n_2}}\\
  &&\rightarrow E[e^{-\lambda\kappa}e^{-(1-\lambda)\tau}].
  \end{eqnarray*}
Thus, letting $\mathbf{n}\rightarrow\infty$ and then $(k_{n_1},k_{n_2})\rightarrow  \mathbf{\infty}$, we have
 $  J_4=o(1)$.

Now, combining the above discussions and letting $\mathbf{n}\rightarrow\infty$ and then $(k_{n_1},k_{n_2})\rightarrow  \mathbf{\infty}$, we get the desired result.

 \textbf{Proof of Theorem 2.2}: Let $\eta_{\mathbf{k}}=\mathbbm{1}_{\left\{\bigcap_{\mathbf{i}\leq \mathbf{k}}\left\{X_{\mathbf{i}}\leq u_{\mathbf{k,i}},X_{\mathbf{i}}(\varepsilon)\leq v_{\mathbf{k,i}}\right\}\right\}}-E\left(\mathbbm{1}_{\left\{\bigcap_{\mathbf{i}\leq \mathbf{k}}\left\{X_{\mathbf{i}}\leq u_{\mathbf{k,i}},X_{\mathbf{i}}(\varepsilon)\leq v_{\mathbf{k,i}}\right\}\right\}}\right)$.\\
Write
\begin{eqnarray*}
&&Var\left(\frac{1}{\log n_1\log n_2}\sum_{\mathbf{k}\in \mathbf{R_n}}\frac{1}{k_1k_2}\eta_{\mathbf{k}}\right)=E\left(\frac{1}{\log n_1\log n_2}\sum_{\mathbf{k}\in \mathbf{R_n}}
\frac{\eta_{\mathbf{k}}}{k_1k_2}\right)^{2}\\
&&=\frac{1}{\left(\log n_1 \log n_2\right)^2}\sum_{\mathbf{k}\in \mathbf{R_n}}\frac{E(\eta_{\mathbf{k}}^2)}{k_1^2k_2^2}+\frac{1}{\left(\log n_1 \log n_2\right)^2}\sum_{\substack{\mathbf{k,l}\in \mathbf{R_n}\\ \mathbf{k} \neq  \mathbf{l}}}\frac{E(\eta_{\mathbf{k}}\eta_{\mathbf{l}})}{k_1k_2l_1l_2}\\
&&=T_1+T_2
\end{eqnarray*}
Since $\eta_{\mathbf{k}}\leq 1$, it follows that
\begin{eqnarray*}
T_1  \leq \frac{1}{(\log n_1\log n_2)^2}\sum_{\mathbf{k}\in \mathbf{R_n}}\frac{1}{k_1^2k_2^2}
 \ll \frac{1}{(\log n_1\log n_2)^2}
\end{eqnarray*}
For $\mathbf{k}\neq \mathbf{l}$ such that $k_1k_2\leq l_1l_2$, we have
\begin{eqnarray*}
\left|E(\eta_{\mathbf{k}}\eta_{\mathbf{l}})\right|&=&\left|Cov\left(\mathbbm{1}_{\left\{\bigcap_{\mathbf{i}\in \mathbf{R_k}}\left\{X_{\mathbf{i}}\leq u_{\mathbf{k,i}},X_{\mathbf{i}}(\varepsilon)\leq v_{\mathbf{k,i}}\right\}\right\}},\mathbbm{1}_{\left\{\bigcap_{\mathbf{i}\in \mathbf{R_l}}\left\{X_{\mathbf{i}}\leq u_{\mathbf{l,i}},X_{\mathbf{i}}(\varepsilon)\leq v_{\mathbf{l,i}}\right\}\right\}}\right)\right|\\
&\leq& \left|Cov\left(\mathbbm{1}_{\left\{\bigcap_{\mathbf{i}\in \mathbf{R_k}}\left\{X_{\mathbf{i}}\leq u_{\mathbf{k,i}},X_{\mathbf{i}}(\varepsilon)\leq v_{\mathbf{k,i}}\right\}\right\}},\mathbbm{1}_{\left\{\bigcap_{\mathbf{i}\in \mathbf{R_l}}\left\{X_{\mathbf{i}}\leq u_{\mathbf{l,i}},X_{\mathbf{i}}(\varepsilon)\leq v_{\mathbf{l,i}}\right\}\right\}}-\mathbbm{1}_{\bigcap_{\mathbf{i}\in \mathbf{R_l}-\mathbf{R_k}}\left\{X_{\mathbf{i}}\leq u_{\mathbf{l,i}},X_{\mathbf{i}}(\varepsilon)\leq v_{\mathbf{l,i}}\right\}}\right)\right|\\
&&+\left|Cov\left(\mathbbm{1}_{\left\{\bigcap_{\mathbf{i}\in \mathbf{R_k}}\left\{X_{\mathbf{i}}\leq u_{\mathbf{k,i}},X_{\mathbf{i}}(\varepsilon)\leq v_{\mathbf{k,i}}\right\}\right\}},\mathbbm{1}_{\bigcap_{\mathbf{i}\in \mathbf{R_l}-\mathbf{R_k}}\left\{X_{\mathbf{i}}\leq u_{\mathbf{l,i}},X_{\mathbf{i}}(\varepsilon)\leq v_{\mathbf{l,i}}\right\}}\right)\right|\\
&\ll& E\left|\mathbbm{1}_{\left\{\bigcap_{\mathbf{i}\in \mathbf{R_l}}\left\{X_{\mathbf{i}}\leq u_{\mathbf{l,i}},X_{\mathbf{i}}(\varepsilon)\leq v_{\mathbf{l,i}}\right\}\right\}}-\mathbbm{1}_{\left\{\bigcap_{\mathbf{i}\in \mathbf{R_l}-\mathbf{R_k}}\left\{X_{\mathbf{i}}\leq u_{\mathbf{l,i}},X_{\mathbf{i}}(\varepsilon)\leq v_{\mathbf{l,i}}\right\}\right\}}\right|\\
&&+\left|Cov\left(\mathbbm{1}_{\left\{\bigcap_{\mathbf{i}\in \mathbf{R_k}}\left\{X_{\mathbf{i}}\leq u_{\mathbf{k,i}},X_{\mathbf{i}}(\varepsilon)\leq v_{\mathbf{k,i}}\right\}\right\}},\mathbbm{1}_{\bigcap_{\mathbf{i}\in \mathbf{R_l}-\mathbf{R_k}}\left\{X_{\mathbf{i}}\leq u_{\mathbf{l,i}},X_{\mathbf{i}}(\varepsilon)\leq v_{\mathbf{l,i}}\right\}}\right)\right|.
\end{eqnarray*}
By the Lemmas 4.2 and 4.1, we get
\begin{eqnarray*}
E\left|\mathbbm{1}_{\left\{\bigcap_{\mathbf{i}\in \mathbf{R_l}}\left\{X_{\mathbf{i}}\leq u_{\mathbf{l,i}},X_{\mathbf{i}}(\varepsilon)\leq v_{\mathbf{l,i}}\right\}\right\}}-\mathbbm{1}_{\left\{\bigcap_{\mathbf{i}\in \mathbf{R_l}-\mathbf{R_k}}\left\{X_{\mathbf{i}}\leq u_{\mathbf{l,i}}, X_{\mathbf{i}}(\varepsilon)\leq v_{\mathbf{l,i}}\right\}\right\}}\right|\ll \frac{l_1l_2-\#(\mathbf{R_l}-\mathbf{R_k})}{l_1l_2}
\end{eqnarray*}
and
\begin{eqnarray*}
&&\bigg|Cov \bigg(\mathbf{\mathbbm{1}}_{\{\bigcap_{\mathbf{i} \in \mathbf{R_{k}}}\{X_{\mathbf{i}}\leq u_{\mathbf{k,i}},X_{\mathbf{i}}(\varepsilon)\leq v_{\mathbf{k,i}}\}\}},\mathbf{\mathbbm{1}}_{\{\bigcap_{\mathbf{i} \in \mathbf{R_{l}}-\mathbf{R_{k}}}\{X_{\mathbf{i}}\leq u_{\mathbf{l,i}},X_{\mathbf{i}}(\varepsilon)\leq v_{\mathbf{l,i}}\}\}}\bigg)\bigg|\\
&&\ \ \ll\left\{
		\begin{array}{lcl}
			\alpha_{\mathbf{l,k},m_{l_1},m_{l_2}}^{*}+
\frac{k_2m_{l_1}+k_1m_{l_2}+m_{l_1}m_{l_2}}{l_1l_2},     &     & {k_1  <     k_2,l_1<l_2};\\
	& &\\		\alpha_{\mathbf{l,k},m_{l_1},m_{l_2}}^{*}+\frac{l_2m_{l_1}+m_{l_1}m_{l_2}}{l_1l_2},   &     & {k_1<l_1,l_2<k_2,k_1k_2<l_1l_2};\\
	& &\\		\alpha_{\mathbf{l,k},m_{l_1},m_{l_2}}^{*}+\frac{l_1m_{l_2}+m_{l_1}m_{l_2}}{l_1l_2},   &     & {l_1<k_2,k_2<l_2,k_1k_2<l_1l_2}.
			\end{array} \right.
	\end{eqnarray*}
respectively.\\
If $k_1<l_1,k_2<l_2$, we have
$$
E\left|\eta_{\mathbf{k}}\eta_{\mathbf{l}}\right|\leq \frac{l_1l_2-\#(\mathbf{R_l}-\mathbf{R_k})}{l_1l_2}+
\alpha_{\mathbf{l,k},m_{l_1},m_{l_2}}^{*}+
\frac{k_2m_{l_1}+k_1m_{l_2}+m_{l_1}m_{l_2}}{l_1l_2}
$$
and
\begin{eqnarray*}
T_2 & \ll& \frac{1}{\left(\log n_1 \log n_2\right)^2}\sum_{\substack{\mathbf{k,l}\in \mathbf{R_n}\\ \mathbf{k} \neq  \mathbf{l}}}\frac{l_1l_2-\#(\mathbf{R_l}-\mathbf{R_k})}{k_1k_2l_1^2l_2^2}
+\frac{1}{\left(\log n_1 \log n_2\right)^2}\sum_{\substack{\mathbf{k,l}\in \mathbf{R_n}\\ \mathbf{k} \neq  \mathbf{l}}}\frac{\alpha_{\mathbf{l},m_{l_1},m_{l_2}}^{*}}{k_1k_2l_1l_2}\\
&&+\frac{1}{\left(\log n_1 \log n_2\right)^2}\sum_{\substack{\mathbf{k,l}\in \mathbf{R_n}\\ \mathbf{k} \neq  \mathbf{l}}}\frac{k_1m_{l_2}+k_2m_{l_1}+m_{l_1}m_{l_2}}{k_1k_2l_1^2l_2^2}\\
&=&T_{21}+T_{22}+T_{23},
\end{eqnarray*}
where
\begin{eqnarray*}
T_{21}&=&\frac{1}{\left(\log n_1 \log n_2\right)^2}\sum_{\substack{\mathbf{k,l}\in \mathbf{R_n}\\ \mathbf{k} \neq  \mathbf{l}}}\frac{l_1l_2-\#(\mathbf{R_l}-\mathbf{R_k})}{k_1k_2l_1^2l_2^2}\\
&=&\frac{1}{\left(\log n_1 \log n_2\right)^2}\sum_{1\leq k_1 \leq l_1\leq n_1}\sum_{1 \leq k_2 \leq l_2 \leq n_2}\frac{k_1k_2}{k_1k_2l_1^2l_2^2}\\
&=&\frac{1}{\left(\log n_1 \log n_2\right)^2}\sum_{l_1=1}^{n_1}\sum_{k_1=1}^{l_1}\sum_{l_2=1}^{n_2}\sum_{k_2=1}^{l_2}\frac{1}{l_1^2l_2^2}\\
&=&\frac{1}{\left(\log n_1 \log n_2\right)^2}\sum_{l_1=1}^{n_1}\sum_{l_2=1}^{n_2}\frac{1}{l_1l_2}
\ll \frac{1}{\log n_1 \log n_2},
\end{eqnarray*}
\begin{eqnarray*}
T_{22}&=&\frac{1}{\left(\log n_1 \log n_2\right)^2}\sum_{\substack{\mathbf{k,l}\in \mathbf{R_n} \mathbf{k} \neq  \mathbf{l}}}\frac{\alpha_{\mathbf{l,k},m_{l_1},m_{l_2}}^{*}}{k_1k_2l_1l_2}\\
&\ll& \frac{1}{\left(\log n_1 \log n_2\right)^2}\sum_{l_1=1}^{n_1}\sum_{k_1=1}^{l_1}\sum_{l_2=1}^{n_2}\sum_{k_2=1}^{l_2}
\frac{1}{k_1k_2l_1l_2\left(\log \log l_1 \log \log l_2\right)^{1+\epsilon}}\\
&=&\frac{1}{\left(\log n_1 \log n_2\right)^2}\sum_{l_1=1}^{n_1}\sum_{l_2=1}^{n_2}\frac{\log l_1 \log l_2}{l_1 l_2 \left(\log \log l_1\right)^{1+\epsilon} \left(\log \log l_2\right)^{1+\epsilon}}\\
&\ll& \frac{1}{\left(\log \log n_1 \log \log n_2\right)^{1+\epsilon}}
\end{eqnarray*}
and
\begin{eqnarray*}
T_{23}&=&\frac{1}{\left(\log n_1 \log n_2\right)^2}\sum_{\substack{\mathbf{k,l}\in \mathbf{R_n}\\ \mathbf{k} \neq  \mathbf{l}}}\frac{k_1m_{l_2}+k_2m_{l_1}+m_{l_1}m_{l_2}}{k_1k_2l_1^2l_2^2}\\
&=&\frac{1}{\left(\log n_1 \log n_2\right)^2}\sum_{l_1=1}^{n_1}\sum_{k_1=1}^{l_1}\sum_{l_2=1}^{n_2}\sum_{k_2=1}
^{l_2}\frac{k_1\log l_2+k_2\log l_1+\log l_1\log l_2}{k_1k_2l_1^2l_2^2}\\
&=&\frac{1}{\left(\log n_1 \log n_2\right)^2}\sum_{l_1=1}^{n_1}\sum_{k_1=1}^{l_1}\sum_{l_2=1}^{n_2}\sum_{k_2=1}
^{l_2}\frac{\log l_2}{k_2 l_1^2 l_2^2}+\frac{1}{\left(\log n_1 \log n_2\right)^2}\sum_{l_1=1}^{n_1}\sum_{k_1=1}^{l_1}\sum_{l_2=1}^{n_2}\sum_{k_2=1}
^{l_2}\frac{\log l_1}{k_1 l_1^2 l_2^2}\\
&&+\frac{1}{\left(\log n_1 \log n_2\right)^2}\sum_{l_1=1}^{n_1}\sum_{k_1=1}^{l_1}\sum_{l_2=1}^{n_2}\sum_{k_2=1}
^{l_2}\frac{\log l_1\log l_2}{k_1k_2 l_1^2 l_2^2}\\
& \ll&\frac{1}{\left(\log n_1 \log n_2\right)^2}\left(\log n_1+\log n_2\right)=o\left(\frac{1}{\left(\log \log n_1\log \log n_2\right)^{1+\epsilon}}\right).
\end{eqnarray*}
Thus
\begin{eqnarray*}
T_2\ll\frac{1}{\left(\log \log n_1\log \log n_2\right)^{1+\epsilon}}.
\end{eqnarray*}
If $k_1<l_1,l_2<k_2$ and $k_1k_2<l_1l_2$
\begin{eqnarray*}
T_2&\ll& \frac{1}{\left(\log n_1 \log n_2\right)^2}\sum_{\substack{\mathbf{k,l}\in \mathbf{R_n}\\ \mathbf{k} \neq  \mathbf{l}}}\frac{l_1l_2-\#(\mathbf{R_l}-\mathbf{R_k})}{k_1k_2l_1^2l_2^2}
+\frac{1}{\left(\log n_1 \log n_2\right)^2}\sum_{\substack{\mathbf{k,l}\in \mathbf{R_n}\\ \mathbf{k} \neq  \mathbf{l}}}\frac{\alpha_{\mathbf{l,k},m_{l_1},m_{l_2}}^{*}}{k_1k_2l_1l_2}\\
&&+\frac{1}{\left(\log n_1 \log n_2\right)^2}\sum_{\substack{\mathbf{k,l}\in \mathbf{R_n}\\ \mathbf{k} \neq  \mathbf{l}}}\frac{m_{l_1}l_2+m_{l_1}m_{l_2}}{k_1k_2l_1^2l_2^2}\\
&=&T_{24}+T_{25}+T_{26}.
\end{eqnarray*}
For the first and second term, we have
\begin{eqnarray*}
T_{24}&=&\frac{1}{\left(\log n_1 \log n_2\right)^2}\sum_{\substack{\mathbf{k,l}\in \mathbf{R_n}\\ \mathbf{k} \neq  \mathbf{l}}}\frac{l_1l_2-(l_1l_2-k_1l_2)}{k_1k_2l_1^2l_2^2}\\
&=&\frac{1}{\left(\log n_1 \log n_2\right)^2}\sum_{l_1=1}^{n_1}\sum_{k_1=1}^{l_1l_2/k_2}\sum_{k_2=1}^{n_2}\sum_{l_2=1}^
{k_2}\frac{1}{k_2l_1^2l_2}\\
&=& \frac{1}{\log n_1\log n_2}=o\left(\frac{1}{\left(\log \log n_1 \log \log n_2\right)^{1+\epsilon}}\right)
\end{eqnarray*}
and
\begin{eqnarray*}
T_{25}&=&\frac{1}{\left(\log n_1 \log n_2\right)^2}\sum_{l_1=1}^{n_1}\sum_{k_1=1}^{l_1l_2/k_2}\sum_{k_2=1}^{n_2}
\sum_{l_2=1}^{k_2}\frac{\alpha^{*}_{\mathbf{l,k},m_{l_1},m_{l_2}}}{k_1k_2l_1l_2}\\
&\ll& \frac{1}{\left(\log n_1 \log n_2\right)^2}\sum_{l_1=1}^{n_1}\sum_{k_1=1}^{l_1l_2/k_2}\sum_{k_2=1}^{n_2}
\sum_{l_2=1}^{k_2}\frac{1}{k_1k_2l_1l_2\left(\log \log l_1 \log \log l_2\right)^{1+\epsilon}}\\
& \ll&\frac{1}{\left(\log n_1 \log n_2\right)^2}\sum_{l_1=1}^{n_1}\sum_{k_2=1}^{n_2}\sum_{l_2=1}^{k_2}\frac{\log l_1+\log l_2-\log k_2}{k_2l_1l_2\left(\log \log l_1 \log \log l_2\right)^{1+\epsilon}}\\
& <& \frac{1}{\left(\log n_1 \log n_2\right)^2}\sum_{l_1=1}^{n_1}\sum_{l_2=1}^{n_2}\frac{\log n_1 \log n_2}{l_1l_2 \left(\log \log l_1 \log \log l_2\right)^{1+\epsilon}}\\
&\leq& \frac{1}{\left(\log \log n_1 \log \log n_2\right)^{1+\epsilon}}.
\end{eqnarray*}
Noting that $n_1=O(n_2)$, we have
\begin{eqnarray*}
T_{26}&\ll&\frac{1}{\left(\log n_1 \log n_2\right)^2}\sum_{l_1=1}^{n_1}\sum_{k_1=1}^{l_1l_2/k_2}\sum_{k_2=l_2}^{n_2}
\sum_{l_2=1}^{n_2}\frac{m_{l_1}}{k_1k_2l_1^2l_2}\\
&<& \frac{1}{\left(\log n_1 \log n_2\right)^2}\sum_{l_1=1}^{n_1}\sum_{k_2=l_2}^{n_2}\sum_{l_2=1}^{n_2}\frac{\log l_1(\log l_1+\log l_2)}{l_1^2 l_2 k_2}\\
& <&\frac{1}{\left(\log n_1 \log n_2\right)^2}\sum_{l_1=1}^{n_1}\sum_{l_2=1}^{n_2}\frac{\log l_1(\log l_1+\log l_2)\log n_2}{l_1^2l_2}\\
&\ll&\frac{\log n_2}{\left(\log n_1 \log n_2\right)^2}\left(\log n_2+\sum_{l_2=1}^{n_2}\frac{\log l_2}{l_2}\right)\\
& \ll& \frac{1}{\left(\log n_1\right)^2}+\frac{\log n_2}{(\log n_1)^2}=o\left(\frac{1}{\left(\log \log n_1 \log \log n_2\right)^{1+\epsilon}}\right).
\end{eqnarray*}
Thus
\begin{eqnarray*}
T_2 \ll \frac{1}{\left(\log \log n_1\log \log n_2\right)^{1+\epsilon}}.
\end{eqnarray*}
If $l_1<k_1,k_2<l_2$ and $k_1k_2<l_1l_2$, by the same arguments as for the second case, we also have
\begin{eqnarray*}
T_2 \ll \frac{1}{\left(\log \log n_1\log \log n_2\right)^{1+\epsilon}}.
\end{eqnarray*}
Therefore
\begin{eqnarray*}
Var\left(\frac{1}{\log n_1\log n_2}\sum_{\mathbf{k}\in \mathbf{R_n}}\frac{1}{k_1k_2}\eta_{\mathbf{k}}\right)\ll \frac{1}{\left(\log\log n_1\log\log n_2\right)^{1+\epsilon}}.
\end{eqnarray*}
Now, the result follows from Theorem 2.1 and Lemma 4.3.

 \textbf{Proof of Theorem 3.1}: To prove the theorem, it suffices to check that $\mathbf{D}(u_{\mathbf{n,i}},v_{\mathbf{n,i}})$, $\mathbf{D}^{*}(u_{\mathbf{k,j}},v_{\mathbf{k,i}}, u_{\mathbf{n,j}},v_{\mathbf{n,i}})$ and
$\mathbf{D}'(v_{\mathbf{n,i}})$  hold.
Recall that $v_{\mathbf{n,j}}\leq u_{\mathbf{n,j}}$ for all $\mathbf{i\leq n}$.
By Normal Comparison Lemma (see e.g., \cite{ Leadbetter_Lindgren_Root1983}),
we have
\begin{eqnarray*}
k_{n_{1}}k_{n_{2}}\alpha_{\mathbf{n},m_{n_1},m_{n_2}}\ll \sum_{\mathbf{1\leq i\neq j\leq n}}|r_{\mathbf{i,j}}|\exp\left(-\frac{v_{\mathbf{n,i}}^{2}+v_{\mathbf{n,j}}^{2}}{2(1+|r_{\mathbf{i,j}}|)}\right),
\end{eqnarray*}
which tends to $0$ as $\mathbf{n}\rightarrow\infty$, by Lemmas 3.3 and 3.4 of \cite{Tan_Wang2014}. Thus, $\mathbf{D}(u_{\mathbf{n,i}},v_{\mathbf{n,i}})$ holds.
Applying Normal Comparison Lemma again, we have
\begin{eqnarray}
\label{proof1}
\sup_{k_{1}k_{2}< n_{1}n_{2}}\alpha_{\mathbf{n},\mathbf{k},m_{n_1},m_{n_2}}^{*}
&\ll& \sup_{k_{1}k_{2}< n_{1}n_{2}}\sup_{\mathbf{I}_{1}\subseteq \mathbf{R_{k}}, \mathbf{I}_{2}\subseteq\mathbf{R_n}\backslash \mathbf{M_{k,n}}} \sum_{\mathbf{i\in I_{1}, j\in I_{2}}}|r_{\mathbf{i,j}}|\exp\left(-\frac{v_{\mathbf{k,i}}^{2}+v_{\mathbf{n,j}}^{2}}{2(1+|r_{\mathbf{i,j}}|)}\right)\nonumber\\
&\ll& \sup_{k_{1}k_{2}< n_{1}n_{2}} \sum_{\mathbf{i\in \mathbf{R_{k}}, j\in\mathbf{R_n}}\atop \mathbf{i\neq j}}|r_{\mathbf{i,j}}|\exp\left(-\frac{v_{\mathbf{k,i}}^{2}+v_{\mathbf{n,j}}^{2}}{2(1+|r_{\mathbf{i,j}}|)}\right)\nonumber\\
&\ll& \sup_{k_{1}k_{2}< n_{1}n_{2}}k_{1}k_{2}\sum_{\mathbf{0\leq j\leq n}, \mathbf{j\neq0}}|\rho_{\mathbf{j}}|\exp\left(-\frac{v_{\mathbf{k,i}}^{2}+v_{\mathbf{n,j}}^{2}}{2(1+|\rho_{\mathbf{j}}|)}\right).
\end{eqnarray}
Now, by a similar arguments as for the proof of Lemmas 3.3 and 3.4 of \cite{Tan_Wang2014}, we can show the term in (\ref{proof1}) tends to $0$ as $\mathbf{n}\rightarrow\infty$. This proves that $\mathbf{D}^{*}(u_{\mathbf{k,j}},v_{\mathbf{k,i}}, u_{\mathbf{n,j}},v_{\mathbf{n,i}})$ holds.
Note that by Normal Comparison Lemma again
\begin{eqnarray*}
|P\left(X_{\mathbf{i}}>u_{\mathbf{n,i}},X_{\mathbf{j}}>u_{\mathbf{n,j}}\right)-(1-\Phi(u_{\mathbf{n,i}}))(1-\Phi(u_{\mathbf{n,j}}))|
\ll |r_{\mathbf{i,j}}|\exp\left(-\frac{u_{\mathbf{n,i}}^{2}+u_{\mathbf{n,j}}^{2}}{2(1+|r_{\mathbf{i,j}}|)}\right),
\end{eqnarray*}
which combined with (\ref{eqTan1}) implies
\begin{eqnarray}
\label{proof2}
&&k_{n_1}k_{n_2}\sum_{\mathbf{i\neq j}\in \mathbf{I}}P\left(X_{\mathbf{i}}>u_{\mathbf{n,i}},X_{\mathbf{j}}>u_{\mathbf{n,j}}\right)\nonumber\\
&&\ll k_{n_1}k_{n_2}\sum_{\mathbf{i\neq j}\in \mathbf{I}}(1-\Phi(u_{\mathbf{n,i}}))(1-\Phi(u_{\mathbf{n,j}}))+
k_{n_1}k_{n_2}\sum_{\mathbf{i\neq j}\in \mathbf{I}}|r_{\mathbf{i,j}}|\exp\left(-\frac{u_{\mathbf{n,i}}^{2}+u_{\mathbf{n,j}}^{2}}{2(1+|r_{\mathbf{i,j}}|)}\right)\nonumber\\
&&\ll k_{n_1}k_{n_2}\left[\sum_{\mathbf{i}\in \mathbf{I}}(1-\Phi(u_{\mathbf{n,i}}))\right]^{2}+
k_{n_1}k_{n_2}\sum_{\mathbf{i\neq j}\in \mathbf{I}}|r_{\mathbf{i,j}}|\exp\left(-\frac{u_{\mathbf{n,i}}^{2}+u_{\mathbf{n,j}}^{2}}{2(1+|r_{\mathbf{i,j}}|)}\right)\nonumber\\
&&\ll \frac{1}{k_{n_1}k_{n_2}}\left[\sum_{\mathbf{i}\in \mathbf{R_{n}}}(1-\Phi(u_{\mathbf{n,i}}))\right]^{2}+\sum_{\mathbf{i\neq j}\in \mathbf{R_{n}}}|r_{\mathbf{i,j}}|\exp\left(-\frac{u_{\mathbf{n,i}}^{2}+u_{\mathbf{n,j}}^{2}}{2(1+|r_{\mathbf{i,j}}|)}\right).
\end{eqnarray}
Since $\sum_{\mathbf{i}\in \mathbf{R_{n}}}[1-\Phi(u_{\mathbf{n,i}})]\rightarrow\tau>0$, the first term in (\ref{proof2}) tends to $0$ as $\mathbf{n}\rightarrow\infty$. By the same arguments as those in the proof of Lemmas 3.3 and 3.4 of \cite{Tan_Wang2014}, the second term in (\ref{proof2}) also tends to $0$ as $\mathbf{n}\rightarrow\infty$. Thus, $\mathbf{D}'(v_{\mathbf{n,i}})$  holds, which completes the proof.

\textbf{Proof of Corollary 3.1}: Letting $v_{\mathbf{n,i}}=x/a_{\mathbf{n}}+b_{\mathbf{n}}+m_{\mathbf{n}}^{*}-m_{\mathbf{i}}$ and $u_{\mathbf{n,i}}=y/a_{\mathbf{n}}+b_{\mathbf{n}}+m_{\mathbf{n}}^{*}-m_{\mathbf{i}}$,
we have
\begin{eqnarray*}
\label{proof3}
&&\bigg(a_{\mathbf{n}}\big(M_{\mathbf{n}}(Z(\varepsilon))-b_{\mathbf{n}}-m_{\mathbf{n}}^{*}\big)\leq x, a_{\mathbf{n}}\big(M_{\mathbf{n}}(Z)-b_{\mathbf{n}}-m_{\mathbf{n}}^{*}\big)\leq y\bigg)\nonumber\\
&&=\bigg(\bigcap_{\mathbf{i\in R_{n}}}\{Y_{\mathbf{i}}(\varepsilon)\leq v_{\mathbf{n,i}}\}, \bigcap_{\mathbf{i\in R_{n}}}\{Y_{\mathbf{i}}\leq u_{\mathbf{n,i}}\}\bigg).
\end{eqnarray*}
Thus, to prove Corollary 3.1, it is sufficient to show the conditions of Theorem 3.1 hold. More precisely, we only need to show that
$\sup_{\mathbf{n}\geq\mathbf{1}}\{n_{1}n_{2}(1-\Phi(v_{\mathbf{n},\mathbf{i}}), \mathbf{i}\leq \mathbf{n})\}$ is bounded, $\sum_{\mathbf{i}\in \mathbf{R_{n}}}[1-\Phi(u_{\mathbf{n,i}})]\rightarrow\tau>0$, $\sum_{\mathbf{i}\in \mathbf{R_{n}}}[1-\Phi(v_{\mathbf{n,i}})]\rightarrow\kappa>0$ and
$\Phi(u_{\mathbf{n,i}})=O(\Phi(u_{\mathbf{n,j}})), \Phi(v_{\mathbf{n,i}})=O(\Phi(v_{\mathbf{n,j}}))$ for $\mathbf{1\leq i\neq j\leq n}$.
It has been done in the proofs of Corollaries 2.2 and 2.3 in \cite{Tan_Wang2014}. The proof is complete.

\textbf{Proof of Theorem 3.2}: The proof is similar with that of Theorem 3.1 by using a comparison lemma for $chi$-random variables to replace the Normal Comparison Lemma, see e.g., \cite{Song_Shao_Tan2022}.

\textbf{Proof of Theorem 3.3}: The proof is similar with that of Theorem 3.1 by using a comparison lemma for Gaussian order statistics to replace Normal Comparison Lemma, see e.g., \cite{Song_Shao_Tan2022}.

\end{document}